\documentstyle[11pt,leqno]{article} 
\input amssym.def
\input amssym.tex

\newtheorem{corr}{Corollary}[section]
\newtheorem{prop}{Proposition}[section]
\newtheorem{theo}{Theorem}[section]
\newtheorem{lemma}{Lemma}[section]
\newtheorem{defi}{Definition}[section]

\newtheorem{rem}{Remark}[section]

\newtheorem{con}{Conjecture}[section]

\newcommand{\spec}{{\rm Spec}}

\begin{document}
\addtolength{\baselineskip}{4pt}

\title{Complex projective surfaces and infinite groups}
\author{F. Bogomolov \thanks{Partially supported by DMS
Grant- 9500774}  and   L. Katzarkov \thanks{Partially
supported by A.P. Sloan Dissertational Fellowship}}
\maketitle

\tableofcontents

\section{Introduction}

It is well known that complex projective surfaces can have highly
nontrivial fundamental groups. It is also known that not all
finitely presented groups can occur as fundamental groups of
projective surfaces. The fundamental problem in the theory then, is to 
determine which groups can occur as the fundamental groups of the complex 
projective surfaces and to describe complex manifolds which occur as 
nonramified coverings
of the surface. The only interesting  case is when the group is infinite 
since finite coverings of  projective surfaces are projective and every  
finite group occurs as a fundamental group of some projective surface.

In this paper we mainly consider a surface with a  representation
as a family of projective curves over a curve. It does not put much
 restriction on the choice of a surface since any surface has such
 a representation after blowing up a finite number of points. We use a base 
change construction combined with a finite ramified covering
to associate to a given surface a collection  of surfaces with infinite
fundamental groups. Most of these fundamental groups were  not  previously 
known to be  fundamental groups of smooth projective surfaces. Every surface
we construct comes equipped with a regular map to a curve of high genus. The 
kernel of the corresponding map of fundamental groups is obtained
from the fundamental group of a generic fiber by imposing torsion relations
on some elements and a very big class of infinite groups occur in this way. 
We also analyze the universal coverings of these surfaces in the context
of the Shafarevich's conjecture which states that the universal covering
of the smooth complex projective variety must be holomorphically convex.

Partial affirmative results \cite{LM}  regarding Shafarevich's conjecture
can be applied to the fundamental groups we obtain. It leads to nontrivial
purely algebraic results on the structure of this groups. On the other hand 
the generality of our construction indicates that not all fundamental groups
we construct will satisfy the algebraic restrictions imposed by Shafarevich's 
conjecture. It would appear then that the
 conjecture may be false in general. We describe a series of potential  
counterexamples in section 4.

We begin with a local version of the general construction.
Namely we have a local fibration without multiple components
 and with only double singular points in the central     fiber. As a first
 step we make a local base change which
 produces a new surface with singular points corresponding to the singular
 points of the central fiber. This move changes the
image of the fundamental group of a generic fiber in the fundamental group
 of the Zariski open subset of nonsingular points. Namely the kernel of this 
map is generated  by the $N$-th powers of the initial vanishing cycle  where 
$N$ is the degree of the local base change. As a second step we desingularize
 the surface by  taking a finite fiberwise covering of the singular
 surface which is ramified at singular points only.

The global construction follows the same pattern, but as a result
 we obtain a surface with a highly nontrivial fundamental
 group coming from the fiber even if the
surface at the beginning was simply connected.
First we construct a singular projective surface with a
big fundamental group of the compliment to the set of singular points.
Smooth projective surface is obtained at the second step as a finite
  covering of the singular surface ramified at singular points only. We prove 
a general theorem (Theorem 2.4) 
which establishes a close similarity between the fundamental groups and 
universal coverings for two
classes of surfaces: normal projective surfaces and the surfaces
 obtained from them by deleting a finite number of points. The above 
construction enlarges
 the image of the fundamental group
of the fiber in the fundamental group of the whole surface. The resulting
 group can be described in purely algebraic terms. Let $\pi_g$ be
 a fundamental  group of a
projective curve of genus $g$  (generic fiber of the fibration).
  Consider  a  finite set
of pairs  of $(s_i \in \pi_g, N_i \in \Bbb{N})$
 and a  subgroup $M$ of the automorphisms    of  $\pi_g$.
We assume that all $s_i$  are vanishing cycles. These are
special conjugacy classes in the fundamental group of the curve
 which constitute a finite number of orbits
 under the action of the mapping class  group $Map(g)$( see section 2).
 The orbits $Ms_i^{N_i}$ generate the normal subgroup $\Xi(M,s_i,N_i)$
 of $\pi_g$.

Now we can give an algebraic version of the description of the
corresponding group.

\begin{defi} Define a Burnside type quotient of $\pi_g$ to be the group
$\pi_g/\Xi(M,s_i,N_i)$.
\end{defi}

In the geometric situation $s_i$ are the vanishing cycles
of the initial fibration and $M$ is its monodromy group.
  Geometric Burnside type groups constitute a proper subset among all
 Burnside type groups. In particular not every data $ s_i,N_i,M$ can be
 geometrically realized. The problem which data can appear in geometry
is the most substantial problem of the above construction.
It is clear that we can change $M$ into its own subgroup of finite index
but at the expense of changing the set of elements $s_i$. The following 
lemma shows that we are  free to vary $N_i$.

\begin{lemma} Let $s_i,N_i = 1 ,M$ describe the data of the smooth
 projective family of curves of  a given genus $g$.
 Then $s_i,N_i ,M$  corresponds to the geometric Burnside group
for any choice $N_i$ satisfying the condition: $N_i = N_j$ if the  singular
 points corresponding to the cycles $s_i$ belong to the same singular
fiber of the fibration.

\end{lemma}

The   topological structure of the algebraic family of curves can
be rather arbitrary if we consider its restriction on a
disc inside  the base curve.
 We summarize relevant results of th article in 
 the following proposition.

\begin{prop} Let $\Gamma$ be any finitely
 presented group and $f_g:\pi_g \to \Gamma$ be any surjective map with
$Ker f_g$ generated by a finite subset $VS$ of the conjugacy classes in $\pi_g$. There exists a surjective map $p_h:\pi_h\to \pi_g$
which corresponds to the contraction of  a set of nonintersecting
handles $H_i$ inside the curve of genus $h$ with the following properties :

1) Every element  $s_j\in VS$ has a preimage $s'_j$ in $\pi_h$ which is
 realized by a smooth  cycle.
 
2) There is a smooth
holomorphic family $X$ of curves of genus $h$ with
 simple singularities over
an algebraic curve  which has all cycles $s_j'$ and the
 generating cycles $a_i,b_i$ of the handles $H_i$ as
 vanishing cycles.

3) The family  $X$ above can have as a monodromy group
  any subgroup of finite index
in $Map(h)$ containing  
 the elements corresponding to 
 positive Dehn twists for all $s_j',a_i,b_i$.
\end{prop}

 Thus any topological data $ s_i, M $ can be realized as part
of some geometric data $ s_i, s_i' M'$ where $M'$ is any subgroup
of finite index in $Map(g)$  containing $M$.

 The previous lemma allows to transfer all additional 
 relations  $s_i'$  into their $N$-powers. Hence  we have two
parameter ``approximation''of any topological data by geometric 
data . One of the parameters corresponds to $N$ and converges to infinity . Another parameter runs through  subgroups of finite index 
in $Map(h)$ which containing the monodromy group $M$ ot the
 fibration over the disc. Since $Map(h)$ is residually
finite the corresponding sequence of groups  converges to $M$.

\begin{rem} These results provide 
with a tool to produce  a big  variety of  fundamental groups
and universal coverings  using  our 
construction. It is worth  noticing that though the initial
group with $N_i = 1$ can be small (even trivial) the  groups
which appear for other choices of $N_i$ are quite diverse.
\end{rem}

\begin{defi} Let us take a free group $\Bbb{F}^{l}$ on
$l$ generators,  the normal subgroup $\Xi$ generated by
 $N$-th powers of  all the primitive
 elements can be described as in definition 1.1 .
 Namely  $\Xi$ is  generated by $Ms^N$ where $s$
is a primitive element  of $\Bbb{F}^{l}$ and $M$
is the group of all automorphisms
of $\Bbb{F}^{l}$. We will denote $\Bbb{F}^l/ \Xi $ 
 by $BT(l,N) $.
\end{defi}

\begin{defi}
Denote by $\pi_{g}/(x^{N}=1)$ the quotient of the fundamental
 group of a Riemann surface of genus $g$ by the group generated by
 the $N$-th powers of all  primitive elements $x$  in $\pi_{g}$. 

\end{defi}

\begin{con}
 (Zelmanov) For big $l ,g,N$  the groups $BT(l,N) $
 and  $\pi_{g}/(x^{N}=1)$ are nonresidually finite groups.

\end{con}

If the above  conjecture is true then the base change   
construction  provides us with many new simple examples of   nonresidually finite 
fundamental groups of smooth projective surfaces.
   The first example of such a group was constructed
 by Toledo (  \cite{TOL}, see also \cite{FL}).

  In section four we discuss the holomorphic convexity of the universal
coverings of the surfaces we have constructed. Recently many new powerful  
methods have been developed to investigate the
 structure of the fundamental groups and universal coverings
 of complex projective surfaces. These methods lead to
many new remarkable results. In particular quite a
 few  positive results on  Shafarevich's conjecture  were
proved e.g.  the results   of
F. Campana,  H. Grauert, R. Gurjar, J. Koll\'ar, B. Lasell, R. Narasimhan,
T. Napier, M. Nori, M. Ramachandran, C. Simpson, S. Shasrty, K. Zuo,
S.T. Yau (see  e.g. \cite{CM}, \cite{LM}, \cite{LN}, \cite{KR}, \cite{K1},
\cite{K2},  \cite{N1}, \cite{SIM}, \cite{YAU}). 

In particular  holomorphic convexity is established (see \cite{LM}) for 
coverings of normal projective surface $X$ corresponding to the homomorphisms
of  $\pi_1(X)$ to $GL(n,\Bbb{C})$ such that the image of these homomorphisms is virtually not equal to $\Bbb{Z}$ (see also \cite{FT}). We  suggest possible 
applications  of the construction related to  the Shafarevich's
 conjecture. Holomorphic convexity of a variety implies
 the absence of infinite chains of compact curves. This can be expressed in 
our case as a restriction on the images on the fundamental groups of
the components of a singular fiber ( see lemma 4.1).  
 We provide with a scheme 
how  to control the behavior of the monodromy and vanishing
cycles to construct a counterexample to the Shafarevich
conjecture.

A free group with a given  number of generators can be identified
with a fundamental group 
of a Riemann surface with one or two ends.
The number of generators in the free group defines the genus of the surface and the number of ends (two if the number of generators
is odd and one if the latter  is even).

\begin{defi}
 Denote by $\Bbb{F}^{g}_k$ a free group with $k$ generators
realized as a fundamental group of a Riemann surface $B$ with one or two ends. 
 Denote by $P^g(k,N)$ the quotient of $\Bbb{F}^{g}_k$ by a normal
 subgroup generated in $\Bbb{F}^{g}_k$ by all primitive elements in
 $\Bbb{F}^{g}_k$ which map into primitive elements (embedded curves) in
$B$.

\end{defi}

 The following group theoretic question is  closely related to the
 Shafarevich conjecture:

{\bf Question} Are there such a $k>1$ and $N$ that that $P^g( 2k,N)$ is a
 finite group and $\pi_{ 2k }/(x^{N}=1)$ is an infinite group?

If the answer to this  question is affirmative then one gets a
  counterexample to the Shafarevich conjecture (see section 4). 

We also suggest potential counterexmaples by considering 
 simplest nontrivial case $N = 3$ in which 
we establish the finiteness of all the groups $P^g(k,3)$( Appendix B). 
In this case we suggest a family of potential counterexamples to the
  conjecture which depends on a non-finiteness of some  groups
in the family ``approximating''  some  infinite group.
 The corresponding family of groups can be described as
 follows. Let us take a chain or ring of curves $X_0$ containing
 more than two curves of genus greater than zero. A natural
contraction of a generic curve of genus $g$ onto this special
fiber defines a subgroup $M_D$ inside the mapping class group
$Map(g)$ which commutes with the contraction.
 If every component $C_i$ of $X_0$ contains at least one 
 vanishing cycle $s_i$ then the group $\pi_g/ ( M_D s_i^3 )$
is infinite though the images of the groups corresponding to different components $C_i\subset X_0$  are finite( Appendix B ). It is true even if we consider smaller subgroups
$M^f \subset M_D$.
Now by our ``approximation'' results we can realize $X_0$ as
a fiber in an  algebraic fibration with the corresponding 
 fiber group $\pi_g / (M^f_j s_i^3 , M^f_j s_{k(N})^{3^N}$. Here
$M^f_j$ is any open subgroup of finite index in $Map(g)$ containing $M^f$,
$N$ is any integer and $s_{k(N)}$ runs through the set
 of additional vanishing classes.
 We conjecture that this set of groups contains many  infinite
groups which will imply that there are many counterexamples to the
 Shafarevich's conjecture obtained from  the surfaces above.

\begin{rem} Since the group we ``approximate'' surjects onto
a nontrivial free product of finite groups we expect that the
 proof of the above conjecture can be found within the 
 modern  theory of infinite groups
\end{rem}

In this section we briefly discuss the symplectic version of our
 construction.  We also formulate an arithmetic variant of the
 Shafarevich conjecture which is presumably easier to prove though the
 conjecture is formally stronger then a direct analogue of the
 complex case.

\bigskip

\noindent
{\bf Acknowledgments:}

The authors would like to thank M.Gromov, J.Koll\'ar, T.Pantev,
G.Tian, D.Toledo for useful conversations and comments.  The second
author would like to thank A.Beilinson, J.Carlson, H.Clemens,
K.Corlette, P.Deligne, R.Donagi, S.Gersten, S.Ivanov, M.Kapovich,
M.Newman, 
\newline M.Nori, A.Olshanskii, M.Ramachandran, C.Simpson, Y.T.Siu,
S.Weinberger, E.Zelmanov and S.T.Yau for  useful conversations and
 constant attention to work.  We would like also to thank the
referee for pointing out an erroneous statement in the initial version
and for helpful suggestions on the organization of the paper. We
 thank 
M.Fried for looking  through the arithmetic part of the
paper.

\section{The general construction}

\subsection{Vanishing cycles - a local construction}

In this subsection we explain the local computation with the vanishing
 cycles on which the whole construction is based. We begin with some
 classical results on degeneration of curves that can be found in
 \cite{DK}.

Let $X_D$ be a smooth complex surface fibered over a disc $D$.  We
 assume that fibers over a punctured disc $D^* = D - 0$ are smooth
 curves of genus $g$ and the projection $ t: X\to D$ is a complex
 Morse function.  In particular the fiber $X_0$ over $0\in D$ has only
 quadratic singular points and it has no multiple components.  Denote
 by $P$ the set of singular points of $X_0$ and by $T : X_t \to X_t$
 the monodromy transformation acting on the fundamental group of the
 general fiber $ X_t $.  This action can be described in terms of Dehn
 twists.  Obviously this action defines an action on the first
 homology group of the general fiber.  The following proposition
 describes completely the topology of $X_D$ and the projection $t: X_D
 \to D$.

\begin{prop}

1) There is a natural topological contraction $ cr: X_D \to X_0$.

2) The restriction of $cr$ to $X_t$ is an isomorphism outside singular
 points $P_i \in X_0$. It contracts the circle $S_i \in X_t$ into
 $P_i$. The monodromy transformation $T$ is the identity outside small
 band $B_i$ around $S_i$ and in $B_i$ the transformation $T$ coincides
 with a standard Dehn twist.

\end{prop}

{\bf Proof} See \cite{DK}.

The above contraction is an isomorphism from $X_t $ minus preimage of
$P$ on $X_{0} - P$.  The preimage of any singular point $P_i$ is a
smooth, homotopically nontrivial curve $S_i\subset X_t$.

\begin{defi} We will call the free homotopy class of $S_{i}$ in 
   the 
 fundamental group of $X_t$ a geometric vanishing cycle.
 It defines a conjugacy class in $\pi_g$(vanishing cycle).
\end{defi}

\begin{rem}

The direction of the standard Dehn twist is defined by the orientation
 of $S_i$ which in turn is defined by the complex structure of the
 neighborhood of the corresponding singular point of the singular
 fiber.
\end{rem}

Let us denote by $De_i$ the topological Dehn transformation of $X_t$
 and by $De_{i,H} $ its action on the homology of $X_{t}$.  In a
 neighborhood of $X_{0}$ the monodromy transformation $T$ is a product
 of Dehn transformations $De_i$ with non-intersecting support.

\begin{lemma}

 1)The monodromy transformation $T$ acts via unipotent transformation
 $T_H$ on the homology group $H_1( X_t, \Bbb{Z})$.

2) $( 1 - T_H)^2 = 0$.

3) $(1 - T_H^N ) = 0(mod N)$ for any $N$.
\end{lemma}

{\bf Proof} It is enough to prove (2) for $De_{i,H}$.  The topological
description of $De_i$ implies that $(1 -De_{i,H})x = (x,s_i)s_i$,
where $s_i$ a homology class of the vanishing cycle $s_i$ and
$(x,s_i)$ is the intersection number.

Since $(s_i,s_i) = 0$ we obtain $(1 - De_{i,H})^2 = 0 $.  The image of
 $(1 - De_{i,H})$ consists of the elements proportional to $s_i$.  We
 also have $(1 - De_{i,H}^N ) = (1 - NDe_{i,H}) = 0 (mod N)$.

The topological transformations $De_i$ commute for different $i$ since
 they have disjoint supports.
  Therefore the same holds for
 $De_{i,H}$.  We also have $ De_{i,H} s_j = 0$ for any vanishing cycle
 $s_j$ since the corresponding circles $S_i, S_j$ don't intersect in
 $X_t$.  We now see that

 \[(1 - T_H)^2 = (1 - \prod De_{i,H})^2 = \prod ( 1 - De_{i,H})^2 = 0
 \]

 where the last equality follows from the above formulas.  Similarly
 we obtain 1) and 3).

$\Box$

The geometric vanishing cycles consist of two different types:

1) The first type includes homologically nontrivial vanishing classes.
 They are all equivalent under the mapping class group $Map(g)$.  The
 latter can be described either as group of connected components of
 the orientation preserving homeomorphisms of the Riemann surface of
 genus $g$ or as the group of exterior automorphisms of the group
 $\pi_g $. The vanishing cycle from this class is a primitive element
 of $\pi_g$ which means it can be included into a set of generators of
 $\pi_g$ satisfying standard relation defining the fundamental group
 of the curve. We shall denote the vanishing cycles of the first type
 type as $NZ$-cycles.

2) The second type consists of elements in $\pi_g$ which are
 homologous to zero. Any vanishing cycle of this type cuts the Riemann
 surface $X_t$ into two pieces and the number of handles in these
 pieces is the only invariant which distinguishes the type of a cycle
 under the action of $Map(g)$. We shall denote the vanishing cycles of
 the second type by $Z$-cycles.

Vanishing $Z$-cycles correspond to the singular points of the singular
 fiber which divide this fiber into two components and $NZ$-cycles to
 the ones which do not.  Each $Z$-cycle defines a primitive element in
 the center of the quotient $\pi_g/[[\pi_g,\pi_g],\pi_g] $, while each
 $NZ$-cycle defines a primitive element in the abelian quotient
 $\pi_g/[\pi_g,\pi_g] = \Bbb{Z}^{2g}$.  Assume now that we made a
 change of the variable $ t = u^N$ and consider the induced family of
 curves over $D$ with a coordinate $u$. Denote the resulting family as
 $X_U$. It is a singular surface and the singular set can be
 identified with the set $P$ o f the singular points of the fiber
 $X_0$. For the new monodromy transformation we have $T^U = T^N$.
 Therefore by lemma 2.3 it acts trivially on $H_1(X_t, \Bbb{Z}_N),
 t\neq 0$.

\begin{lemma} The surface $X_U$ contracts to the central fiber $X_0$.
\end{lemma}

 Indeed the fiberwise contraction of $X_D$ to $X_0$ can be lifted into
 a contraction of $X_U$.

\begin{rem} The fundamental group
 $\pi_1(X_D - P) = \pi_1(X_D) =\pi_1(X_0)$ since the singular points
of the fiber $X_0$ are nonsingular points of $X_D$.  The analogous
statement is not true however for $X_U$.
\end{rem}

\begin{theo} The fundamental group $\pi_1(X_U - P)$ is
equal to the quotient of $\pi_1(X_t) = \pi_g$ by a normal subgroup
 generated by the elements $s_i^N$.
\end{theo}

{\bf Proof} The fundamental group of a complex coincides with the
 fundamental group of any two-skeleton of the complex. There is a
 natural two-dimensional complex with a fundamental group as in the
 theorem. Namely let us take a curve $X_t$ and attach two-dimensional
 disks $D_i$ via the boundary maps $f_i : dD_i\to S_i$ of degree $N$.
 The resulting two dimensional complex $X_t^a$ evidently has a
 fundamental group isomorphic to the group described in the
 theorem. We are going to show that $X_t^a$ can be realized as a two
 skeleton of $X_U - P$. We prove first the following statement:

\begin{lemma} The surface $X_U - P$ retracts onto a three
dimensional complex $X^3$ which is a union of generic curve $X_t$ and
 a set of three-dimensional lens spaces $L_N^i$. Each lens space
 corresponds to the singular point $P_i$ of the singular fiber and
 $L_N^i$ intersects $X_t$ along the band $B_i$.
\end{lemma}

{\bf Proof }Locally near $P_i$ the surface $X_U$ is described by the
equation $ t^N = z_1z_2 $. Hence a neighborhood $U_i$ of $P_i$ is a
cone over the three-dimensional lens space $L_N^i= S^3/
\Bbb{Z}_N$. Here $S^3$ is a three dimensional sphere and $\Bbb{Z}_N$
is generated by the matrix with eigenvalues $\chi, \chi^{-1}$, where
$\chi$ is a primitive root of unity of order $N$. We can assume that
the surface $X_t$ intersects $L_N^i$ along a two-dimensional band
$B_i$ with a central circle $S_i$ defining a generator of
$\pi_1(L_N^i)$. Though the band $B_i$ is a direct product of $S_i$ by
the interval its embedding into $L_N^i$ is nontrivial : the boundary
circles have nonzero linking number. We use a fiberwise contraction to
contract $X_U - P$ to the union to $X^3$. It coincides with a standard
contraction outside of the cones over $L_N^i$. Therefore we obtain a
contraction of $X_U -P$ onto $X^3$.

$\Box$

The two-skeleton of $X^3$ can be obtained as a union of $X_t$ and a
two-skeleton of $L_N^i$. The latter can be seen as a retract of a
complimentary set to the point in $L_N^i$. Here is the topological
picture we are looking at.

\begin{lemma} Let $L_N^i$ be a lens space as  above. Then the complimentary set to a point retracts on  complex $L^2_i$  obtained by attaching a disc to the
circle $S_i$ via the boundary map of degree $N$.
\end{lemma}

{\bf Proof} The sphere $S^3$ can be represented as a joint of two
 circles $S_i$ and $S'$. In other words it consists of intervals
 connecting different points of $S_i, S'$.  The action of $\Bbb{Z}_N$
 with $ L_N^3 = S^3/ \Bbb{Z}_N$ rotates both circle. Define the disc
 $D_x$ to be a cone over $S_i$. Different discs $D_x$, $x \in S'$,
 don't intersect and the image of $D_x$ in $L_N^3$ is $L^2_i$. The
 fundamental domain of $ \Bbb{Z}_N$-action lies between discs $D_x
 ,D_{gx}$, where $g$ is a generator of $\Bbb{Z}_N$. Hence this domain
 is isomorphic to $D^3$ and coincides with complimentary of $L^2_i$ in
 $L_N^i$. The last proves the lemma.

$\Box$

\begin{corr} There is an embedding of $X_t^a $ into $X_U - P$ which induces an isomorphism of the fundamental groups.

\end{corr}

 Indeed we obtain the two skeleton of $X^3$ gluing $X_t$ and $L^2_i$
 along $S_i$, but the resulting two-complex coincides with
 $X_t^a$. Since $X^3$ is a retract of $X_U - P$ we obtain the
 corollary and finish the proof of the theorem.

$\Box$

Let us denote by $G_X$ the fundamental group of $X_U - P$ and by
$\widetilde X$ its universal covering. The description of the group
$G_X$ can be obtained in pure geometric terms. Namely the vanishing
cycles $S_i$ don't intersect and therefore we can first  contract
$Z$-cycles to obtain the union of smooth Riemann surfaces $X_i$ with
normal intersections.  The graph corresponding to this system of
surfaces is tree since every point corresponding to $Z$-cycle splits
it into two components. Remaining $NZ$- cycles lie on different
surfaces $X_j$ and constitute a finite isotropic subset of primitive
elements in $\pi_1(X_j)$. The vanishing cycle $S_i$ which contracts to
the point of $X_j$ defines a map of a free group onto $\pi_1(X_j)$
with $S_i$ corresponding to standard relations in $\pi_1(X_j)$. We
also have the following lemma.

\begin{lemma} If $N$ is odd or divisible by four then the group $G_X$ has a natural surjective projection on the quotient group $\pi_g/[[\pi_g,\pi_g],\pi_g]$ with additional relation $x^N = 1$.
If $N = 2$ then $G_X$ maps surjectively on a central extension of
$\Bbb{Z}_2^{2g}$ by $\Bbb{Z}_2$ and the images of all $s_i$ have order
$2$.
\end{lemma}

{\bf Proof} The case when $N$ is odd or divisible by four is clear
since all the elements $s_i$ have order $N$.  In the case $N = 2$ we
use the fact that $NZ$ cycles $s_i$ contained in the component $X_j$
lie in the isotropic subspace. Hence there exists a standard
$\Bbb{Z}_2$ central extension of $\Bbb{Z}_2^{2g_j}$ where the images
of all $s_i\in X_j$ are exactly of order $2$.  Denote this group by
$G_j$ and the generator of its center by $c_j$. Now consider the
product of $G_j$ for all $j$ and factor it by a central subgroup
generated by $ c_k - c_j$ if $X_k$ and $X_j$ intersect. Since the
graph is a tree we obtain that the quotient of the center by the above
group is equal to $\Bbb{Z}_2$ which is identified with the zero
homology with coefficients $\Bbb{Z}_2$.  The image of $\Bbb{Z}$-cycle
$s_i$ coincide with $c_j$ if $s_i\in X_j$ and hence is never
zero. Vanishing $NZ$-cycles project into nonzero elements of the
abelian quotient of the group.  $\Box$

\begin{rem} If $N$ is odd or divisible by four  we obtain a canonical quotient
of $G_X$ which is equivariant under $Aut(\pi_g)$. We denote this group by
 $UC_g^N$. For $N = 2$ our construction is less canonical,
 it depends on the choice of isotropic subspace in $H_1(X_t,
 \Bbb{Z}_2)$

\end{rem}

Here we start to develop an idea that will be constantly used through
 out the paper. Namely we show that we can work with open surfaces and
 get results concerning the universal coverings of the closed
 surfaces.

\begin{theo} There is a natural $G_X$-invariant
embedding of $\widetilde X$ into a smooth surface $\widetilde X_U $
 with $\widetilde{X}_U/G_X = X_U$. The complimentary of $\widetilde X$
 in $\widetilde X_U$ consists of a discrete subset of points.
\end{theo}

{\bf Proof} Let $P_i $ is a point that corresponds to a $NZ$ vanishing
 cycle. The preimage of a local neighborhood $U_i $ of $P_i$ in
 $\widetilde X$ consists of a number of nonramified coverings of $U_i
 - P_i$.  Since $\pi_1(U_i - P_i) = \pi_1(L_N^i) = \Bbb{Z}_N$ we
 obtain that nonramified covering we get is finite.  Therefore by
 local pointwise completion we obtain $\widetilde X_U$ with the action
 of $G_X$ extending on it.

 Now the universal nonramified covering of $U_i - P_i$ is a punctured
 unit disk in $\Bbb{C}^2$. Therefore if the map $\pi_1(U_i - P_i)\to
 G_X$ is injective then the surface $\widetilde X_U$ is nonsingular.
 Since the group $\pi_1(U_i - P_i) = \Bbb{Z}_N$ and is generated ed by
 the vanishing cycle $s_i$ we can easily get the result.  Indeed the
 groups $H_1(G_X, \Bbb{Z}_N)$ and $H^1( X_t, \Bbb{Z}_N)$ are
 isomorphic.  The general result follows from theorem 2.1 since we
 have constructed finite quotients for $G_X$ which map every $s_i$
 into an element of order $N$.

$\Box$

We also have the following:

\begin{lemma} Let $G_X^0$ be the kernel of projection
 of $G_X$ into one of the finite groups defined in lemma 1.13.  Then
 $G_X^0$ acts freely on $\widetilde X_U$ and the quotient is a family
 of compact curves without multiple fibers and with a fundamental
 group $G_X^0$.
\end{lemma}
 {\bf Proof} Indeed the surface $\widetilde X_U$ was obtained from
 $\widetilde X$ by adding points and since $\widetilde X$ was
 simplyconnected the former is symplyconnected either.  Any element of
 $G_X$ which has invariant point in $\widetilde X_U$ is conjugated to
 the power of $s_i$, but the latter are not contained in $G_X^0$.

$\Box$

We are done with our description of the local computations.  By taking
 finite coverings and taking away the singularities of the covering we
 were able to make all local geometric vanishing cycles to be torsion
 elements.

\subsection{ The global construction}

In this subsection we globalize the construction of section 2.1
to a compact surface.

Let  $X$ be smooth surface with a proper map to a smooth projective
 curve $C$. We assume that the map $f:X\to C$ is described locally by
 a set of holomorphic Morse functions and hence satisfies  the
 conditions of the previous section  in the neighborhood of any fiber.   The generic fiber is a
 smooth curve $X_t, t\in C$ of genus $g$.

 We denote by $P$ the set of all singular
points of the fibers and by $P_C$ the set of points in $C$
corresponding to the singular fibers , $ f(P) = P_C$. The main difference of the
global situation lies in the presence of the global monodromy group
which is the image of $\pi_1( C - P_C)$ in the mapping class group
$Map(g)$.  We denote this group by $M_X$.

Let us choose an integer $N$ and consider a base change $h: R\to
 C$ where $R$ such that the map $h$ is $N$-ramified at all the
 preimages of the points from $P_C$ in $R$. Consider a surface $S$
 obtained via a base change $h: R\to C$. We have the finite map $h' :
 S \to X$ defined via $h$ and the projection $g: S\to R$ with a
 generic fiber $S_t = X_{h(t)} $.  The surface $S$ is singular with
 the set of singular points equal to $h^{-1}(P) = Q$ and the set of
 singular fibers over the points of $h^{-1}P_C = P_R $..  The monodromy
 group $M_S$ of the family $S$ is subgroup of finite index of the
 group $M_X$.

\begin{theo} The fundamental group $\pi_1(S - Q)$ surjects onto
 $\pi_1(R)$. The kernel of this surjection is a  quotient of $\pi_g$ by a normal subgroup
generated by the orbits of $M_S (s_i^N)$.
\end{theo}

{\bf Proof} Indeed we have a natural surjection of $\pi_g$ on the
 kernel of projection onto $\pi_1(R)$ since there are no multiple
 fibers in the projection $g$. A standard 
 argument reduces all relations to the local ones and the local
 relations where described in the previous section (see theorem 2.1).

$\Box$

Now we move to the second step of our construction getting out of $S
 -Q$ a smooth compact surface $S^{N}$ with almost the same fundamental
 group.  In the next subsection we will develop some partial theory of
 this second step summarizing and generalizing some known results. Let
 us assume that $N$ is either odd or divisible by $4$.

\begin{lemma} There exists a smooth projective surface
 $S^{N}$ with a finite map $f : S^{N}\to S$ such that the image of the
 homomorphism $f_{*}:\pi_1(S^{N})\to \pi_1(S -Q)$ is a subgroup of
 finite index in $\pi_1(S -Q)$.
\end{lemma}

 As it was already shown in lemma 2.5 there exists a projection of
 $\pi_1(X_g)$ to a finite group $UC_g^N$ which is invariant under
 $Aut(\pi_g)$ and  factors through the fundamental group of a
 small neighborhood of a singular fiber in $S-Q$.  Therefore, if the
 map $h: S \to R$ has a topological section we obtain a finite
 fiberwise covering of $S$ which is ramified only over the singular
 set $Q$.  The resulting surface coincides locally with the smooth
 surface described in theorem 2. It is also smooth and the map $f$ is
 finite.  If there is not a topological section we first make a base
 change nonramified at $P_R$ in order to obtain such a section ( see
 theorem 5.2 Appendix A) and then apply the previous argument.  Since
 the map $f$ is finite the preimage $f^{-1}(Q)$ consists of a finite
 number of smooth points and therefore $ \pi_1(S^{N}) = \pi_1( S^{N} -
 f^{-1}(Q)) $ (see theorem 5.1 Appendix A). This also proves 
 that the image of the homomorphism $f_{*}: \pi_1(S^{N})\to \pi_1(S
 -Q)$ is a subgroup of finite index in $\pi_1(S -Q)$.

$\Box$

As we have said in the introduction to find a counterexample to the
 Shafarevich's conjecture we try to control the existence of an
 infinite connected chain of compact curves in the universal covering
 of $S^{N}$. The above construction shows that we can do it by
 controlling the image of the fundamental groups of the open
 irreducible components of the reducible singular fibers in the
 fundamental group of the open surface $S-Q$.

\subsection{Comparison theorem}

The construction discussed in the previous section can be applied to
both projective and quasiprojective fibered surfaces. This indicates
that the universal coverings and the fundamental groups for these two
classes of surfaces have a similar structure. In this section we
illustrate another flavor of the same principle by describing a
procedure comparing the fundamental groups of surfaces with quotient
singularities to the fundamental groups of certain smooth surfaces.
For future reference we will set up this transition in a slightly
bigger generality.

Let $V$ be a normal projective complex surface and $Q \subset V$ be
the finite set of its singular points. The fundamental group
$\pi_1(V)$ is the quotient of $\pi_1(V - Q)$ by the normal subgroup
generated by the images of the local fundamental groups of the points
$q\in Q$.  Recall that the local fundamental group $L_{q}$ of a point
$q \in Q$ is defined as the fundamental group of a deleted
neighborhood of $q$, i.e.  as the group $\pi_1( U_q - \{ q \})$ where
$U_q$ is a small analytic neighborhood of the point $q\in V$.

The topology of the neighborhood $U_q$ is completely determined by
$L_q$ as it was shown by D. Mumford \cite{MUM}. The following theorem
is a generalization of a theorem by J. Koll\'ar (see \cite{K2}).

\begin{theo}
Let $V$ be a normal projective surface and let $Q$ be the finite set
of its singular points. Consider for any $q\in Q$ a normal subgroup of
finite index $K_q\triangleleft L_q$ which contains the kernel of the
natural map $ L_q\to \pi_1(V - Q)$. Then there exists a smooth
projective surface $F$ and a surjective finite map $r : F \to V$ which
induces an isomorphism between $\pi_1(F)$ and the quotient of $\pi_1(V
- Q) $ by the normal subgroup generated by the images of $K_q \in
\pi_1(V - Q) , q\in Q$.
\end{theo}
{\bf Proof.}  Denote by $K_Q$ the normal subgroup of $\pi_1(V - Q)$
generated by the subgroups $K_q\subset L_q$.  We obtain the surface
$N$ as a generic hyperplane section of a singular projective variety
$W$ with the following property: $W$ contains a subvariety $S$ of
codimension $\geq 3$ with $\pi_1(W-S) = \pi_{1}(V-Q)/K_{Q}$. We may
assume that $F$ does not intersect $S$ since the latter has
codimension at least $3$ in $W$. The fundamental group $\pi_1(F) =
\pi_1(W-S)$ since $F$ is a generic complete intersection in $W$. We
are going to construct $W$ as a union of two quasiprojective
subvarieties. Denote by $G_q$ the finite quotient $L_q/K_q$ and by $G$
the direct product of all the groups $G_q$. Denote by $g_q$ the
coordinate projection of $G$ onto $G_q$ and by $i_q$ the coordinate
embedding of $G_q$ into $G$. For any $q$ there is a natural finite
covering $ M_q$ of $U_q$ corresponding to the projection $L_q\to
G_q$. The preimage of $q$ in $M_{q}$ consists of a single point and
the projection $M_q\to U_q$ is nonramified outside $q$.  In the next
lemma we prove the existence of an algebraic extension of this local
covering.

 \begin{lemma} There exist an open affine subvariety $V_q\subset V$
 containing $q$ and an affine variety $B_q$ which is a $G_{q}$-Galois
 covering of $V_q$ ramified only at $q$ so that $B_{q}\times_{V_{q}}
 U_{q}$ and $M_{q}$ are isomorphic.
\end{lemma}
{\bf Proof} Let $\hat A_q$ be the completed local ring of $q\in V$.  A
local $G_q$-covering defines a finite algebraic extension $\hat B_q$
of $\hat A_q$. By Artin's approximation theorem \cite{ART} there
exists an affine ring $A\subset {\Bbb C}(V)$ and a finite algebraic
extension $B$ of $A$ which locally at $q$ corresponds to the extension
$\hat B_q$ over $\hat A_q$.  Explicitly the extension $\hat B_q$ is
described by a monic polynomial $f(x)$ with coefficients in $\hat
A_q$. If we now consider any monic polynomial $g(x)$ over the ring $A$
with $g(x) = f(x) \; {\rm mod} \; {\frak m}_q^N $ for a big enough $N$
then the resulting algebraic extension $B$ will be the one we need.
The ring $A$ defines an open algebraic subvariety $\spec(A) \subset V$
containing $q$.  Similarly $B$ defines an affine variety $\spec(B)$
with a finite projection $p_q : \spec(B) \to \spec(A)$. This
projection is unramified outside $q$ in the formal neighborhood of
$q$.  Since we have the freedom to impose any finite number of extra
conditions on $g(x)$ we can choose $p_q $ to be unramified at any
finite number of points. In particular we may assume that the
projection $p_q$ is nonramified over $Q-\{ q\}$. That means that the
divisorial part $D\subset \spec(A)$ of the ramification of $p_q$ does
not intersect $Q$.  Now we can take $\spec(A)/D$ as $V_q$. Let $B_q $
denote $\spec(B) \times_{\spec(A)} V_{q}$.  It is an affine variety
with affine $G_q$-action since it extends the local nonramified Galois
covering $U_q-\{ q\}$ and has the same degree.  \hfill $\Box$

\medskip

Let $B_{0}$ be the product of all $B_{q}$'s over $V$.  This is an
affine variety with the action of $G$. The quotient $B_0/G = V_0'$ is
an open affine subvariety of $V$ which contains $Q$.  The action of
$G$ on $B_0$ is free outside of the preimage of $Q$. Let $G \to GL(E)$
be a (not-necessarily irreducible) faithful linear representation of
$G$ of dimension $e$ with the property that only $1\in G$ is
represented by scalar matrix. Consider the diagonal action of $G$ on
the product $B_0 \times E$. There exists a natural $\Bbb{C}^*$ action
on $E$ - multiplication by scalars.  It extends to a
$\Bbb{C}^*$-action on the product which commutes with the $G$-action.

Let $F_0' = (B_0\times E)/G$ be the quotient variety. It is an affine
variety with induced $\Bbb{C}^*$ action which has a natural projection
$\pi_0 : F_0' \to V_0'$ and a zero section $i(V_0') = (B_0\times
0)/G$.  For any $s\in V_0'-Q$ the preimage $\pi_0^{-1}(s)$ is a vector
space isomorphic to $E$. Moreover $F_0'$ contains a natural vector
bundle $I$ over $V_0'-Q$.  Its sheaf of sections coincides with the
sheaf of $G$-equivariant sections of the constant sheaf ${\cal
O}\otimes E$ over $B_0$.  Let us choose a smaller affine variety
$V_0\subset V_0', Q\subset V_0$ with the property that $I$ is constant
on $V_0-Q$. We define $F_0$ as the preimage of $V_0$ in $F_0'$.  Let
$V_1$ be an open subvariety of $V$ which does not contain $Q$ and such
that the union of $V_1$ and $V_0$ is equal to $V$. Let $J$ be the
trivial bundle of rank $e$ on $V_1$. Choose a linear algebraic
isomorphism of $F_0$ and $J$ over the intersection of $V_0$ and
$V_1$. Use this isomorphism to glue the projectivization ${\Bbb P}(J)$
with the singular variety $X = (N_0 -i(V_0')/{\Bbb C}^*$. The
resulting proper variety $W$ has a natural projection $p: W\to V$ with
all the fibers outside $Q$ isomorphic to projective space
$\Bbb{P}^{e-1}$. Moreover the preimage of $V-Q$ in $W$ coincides with
the projectivization of a vector bundle according to the construction
of $W$. The fiber $W_q$ over $q$ coincides with
$\Bbb{P}^{e-1}/i_q(G_q)$.

The action of $i(G_q)$ on $\Bbb{P}^{e-1}$ is effective because of our
 assumption on the representation of $G \to GL(E)$.  We denote by
 $S_q$ the singular subset of $W_q$.. It lies in the image of the fixed
 sets ${\rm Fix}_{g}(\Bbb{P}^{e-1})\subset \Bbb{P}^{e-1}$ for
 different elements $g\in i_q(G_q), g\neq 1$. Define a subvariety $S$
 as the union of the varieties $S_q , q\in Q$. The set $S$ has
 codimension $\ge 3$ in $W$ since the codimension of $S_q$ in $W_q$ is
 at least $1$.

\begin{lemma}  The variety $W-S$ has a fundamental group isomorphic to 
$\pi_1(V-Q)/ K_Q$.
\end{lemma}
{\bf Proof} The fundamental group of $W-S$ is the quotient of the
 fundamental group of $\pi_1(V-Q)$ since $W$ contains an open
 subvariety which is $\Bbb{P}^{e-1}$ fibration over $V-Q$ and
 therefore has the same fundamental group. The group $K_q$ maps into
 zero under the surjective map $\pi_1(V-Q)\to \pi_1(W-S)$ since the
 image of a neighborhood of $q$ via the zero section $i$ has $K_q$ as
 a local fundamental group.  All the relations are local and
 concentrated near special fibers.  A formal neighborhood of $W_q -
 S_q$ in $W - S$ is topologically isomorphic to a fibration over $W_q
 - S_q$ with $M_q$ as a fiber.  Therefore all local relations follow
 from $K_q =1$.  It finishes the proof that $\pi_1( W - S) = \pi_1(V -
 Q)/K_Q$. \hfill $\Box$

Finally we prove the projectivity of $W$ by constructing an ample line
bundle on it. Start with a line bundle $L$ on $W$ whose sections give
an embedding of $X$ into a projective space. To see that such an $L$
exists consider first the $G$-invariant and $\Bbb{C}^*$ homogeneous
sections of the trivial bundle ${\cal O}\otimes E$ over $B_0$. Assume
that the degree of homogeneity is big enough and divisible by the
order of $G$.  It is known that such sections separate the points in
the quotient variety $X$ and we obtain and embedding of $X$ into a
projective space. Thus the induced bundle ${\cal O}(1)$ is defined on
$X$.  Denote by $L$ some extension of ${\cal O}(1)$ to $W$.  Such
extension exists since the complement of $X$ in $W$ is smooth.

Next by choosing a polarization $H$ on $V$ appropriately we may assume
that the global sections of $L\otimes p^*H$ on $W$ separate all the
points of $X$. The restriction of $L$ on ${\Bbb P}(J)$ coincides with
${\cal O}_{{\Bbb P}(J)}(m)$ for some positive integer $m$. By
replacing $H$ with a high power of $H$ if necessary we can produce
enough sections of $L \otimes p^*H$ to separate the points of ${\Bbb
P}(J)$.  Therefore $L \otimes p^*H$ gives an embedding of $W$ into a
projective space.  \hfill $\Box$

\begin{rem} J.Koll\`{a}r   (\cite{K1})
obtained similar result under additional assumption of existence of a
surjective map of $\pi_1(V - Q)$ onto a finite group $H$ with $K_q$ as
a kernel of the induced map on $L_q$ for every $q\in Q$.  \end{rem}

 The universal covering $\tilde F$ of the smooth projective surface
$F$ is very similar to the ramified covering $\tilde V$ of $V$
corresponding to the quotient group $\pi_1(V - Q)/K_Q$. Namely there
is a natural finite map $p : F \to V$ which induces a finite map from
$\tilde F$ to $\tilde V$. Hence both $\tilde F$ and $\tilde V$ are
simultaneously either holomorphically convex or not.

\subsection{Fiber groups}

 The global construction described in section 2.2
 treats separately the part of
 the fundamental group of the fibered surface which lies in the image
 of the fundamental group of the fiber.  Let $V$ be normal projective
 complex surface and $Q$ be a set of its singular points. Suppose that
 there is a projection of $V$ on a smooth curve $C$ which has no
 multiple fibers and the generic fiber of the projection is a curve of
 genus $g > 1$.

\begin{defi} Denote by $\pi_{1,f} (V - Q) $
  the image of the fundamental group $\pi_g$
of the general 
 fiber in $\pi_1(V - Q)$. We will call this 
 group  a general fiber group.
\end{defi}

In this article we mostly consider the case when the set of singular
 points in $V$ includes only singularities with finite local
 fundamental groups. It is well known that these are exactly the
 quotient singularities.

\begin{defi} We will call  the group $\pi_{1,f}(V - Q)$ above
a fiber group if $Q$ consists of the quotient singularities only.

\end{defi}

 We shall also give a special notation for the case when $Q$ is empty.

\begin{defi} We will call  the group $\pi_{1,f}(V)$ projective fiber
 group if $V$ is a projective fibered surface over $C$ without
 multiple fibers.

\end{defi}

\begin{rem} Though we don't  allow multiple fibers in the above definition
 we allow some multiple components in the singular fibers.  We need
 that at least one component of each singular fiber has multiplicity
 one.

\end{rem}

Thus we have defined three classes of groups. These groups are
 equipped with a surjective map from the group $\pi_g$.  The principal
 difference between these three classes of groups lies in the
 nontriviality of the local fundamental groups of normal surface
 singularity.

\begin{rem}
It follows that the general fiber group occurs also as a fiber group
 of a projective surface if the images of all local fundamental groups
 are finite.  In particular this is true if all singular points have
 finite local fundamental groups.

\end{rem}

As a consequences of theorem 2.4 we get:

\begin{corr}
The classes of fundamental groups and fiber groups are the same for
projective smooth surfaces and projective surfaces minus quotient
singularities.

\end{corr}

The above results suggest that finding examples of smooth projective
 surfaces with pathological behavior of fundamental groups and
 universal coverings can be reduced to a similar problem for normal
 projective surfaces minus singular points. The latter seems to be an
 easier task.

\section {Nonresidually finite groups}

This section contains some material that shows opportunities to make our construction applicable to  a big variety of examples. The second subsection shows how one can use our construction  and a  conjecture by Zelmanov to obtain a variety of potential ex
amples of surfaces with nonresidually finite fundamental groups. It is a pleasure for the second author  to thank M. Nori for many illuminating discussions concerning the related ideas.

\subsection{ Variety  of constructions} 

In this subsection we analyze the groups that can be obtained as fiber
groups.  Recall that the fiber group depends on the genus $g$ of the
generic fiber, the monodromy group and the set of geometric vanishing
cycles.  Thus we can define an abstract algebraic data which gives us
an abstract analogue of the fiber group.  Let $M$ be a subgroup of the
mapping class group $Map(g)$ and $g_1,...,g_N$ be any finite
collection of elements in $\pi_g$ which are powers of vanishing
cycles. Let $Mg_i$ be the orbit of $g_i$ under the action of $M$.

\begin{defi}

An abstract fiber invariant is a set of the form $(g,M, Mg_i)$ for
some $M\subset Map(g)$ and some finite set of $g_i\in \pi_g$ as above.

\end{defi}

An abstract fiber invariant defines an abstract fiber group as the
quotient of $\pi_g$ by a normal subgroup generated by $Mg_i$.  We are
interested in determining conditions under which this abstract fiber
group is the actual fiber group of some geometric fibration which
means a complex quasiprojective surface with a fibration over a curve.

\begin{rem} The answer is rather simple in smooth or symplectic
categories because all the elements of $Map(g)$ can be realized by the
automorphisms of the Riemann surface of genus $g$ which preserve a
given volume form. However the question about geometric fiber
invariants (smooth projective case) is substantially more delicate.
\end{rem}

The following theorem shows that we still have a significant freedom
to vary geometric fiber invariants.

\begin{theo} Assume that $(g, M, Ms_i)$ is a fiber invariant
of a projection $p: X \to C$ where $X$ is a smooth compact surface and
$p$ has only Morse singularities and $s_i$ are the vanishing cycles
corresponding to the singular points $P_i$ of the fibers of $p$. Then
$(g,M,Ms_i^{N_i})$ is also a fiber invariant associated with some
other fibration provided $N_i = N_j$ if $P_i, P_j$ are contained in
the same singular fiber.

\end{theo}

{\bf Proof} We shall construct a new fibration with the desired
 properties by applying the base change construction to the fibration
 $p : X \to C$.  Let $P$ be the set of points of $C$ corresponding to
 the singular fibers of $p$ and $N(p)$ be a function on $P$ with
 positive integer values such that $N(s) = N_i = N(P_i)$ if a singular
 point $P_i\in X_s$ for some $s\in P$.

 Let us take a base change $ h : R\to C$ where $R$ is a cyclic
covering of $C$ of degree $N$ - the  minimal integer divisible by
all $N_i$, and with ramification indices $N_i$ at $P_i$ and $N$ at some
$c\notin P$.

This covering satisfies the following properties:

1) The  preimage of $s\in P$ in $R$ is $N(s)$ ramified.

2) The map $\pi_1(R - h^{-1}(P)) \to \pi_1(C - M) $ is surjective.

Indeed the first property is obvious from the construction of $h$.
 The point $c$ has exactly one preimage in $R$.  Hence any closed loop
 in $C - P$ containing $c$ lifts into a closed loop in $R - h^{-1}(P)$
 which proves surjectivity of the corresponding map of fundamental
 groups.

Now we can induce the family of curves on $R$ by the map $h$.  The
resulting surface $Y$ is a singular surface with a finite map $f :
Y\to X$. The singular points of $Y$ are the preimages of the points
$P_i\in X$ with $N_i > 1$. All the singularities of $Y$ are quotient
singularities. If we denote the set of singular points by $Q$ then the
fiber invariant of the projection $ p_h : ( Y - Q) \to R$ is described
by the data $ g, M, Mg_i$.

Indeed the monodromy of the new family coincides with the image of the
 group $\pi_1 ( R - h^{-1}( P)$ in $Map(g)$ but the latter is equal to
 the image of $\pi_1(C - P)$ in $Map(g)$ as it was proved above.
 Hence the monodromy of the newly obtained family is $M$. All
 relations in the fiber group of $Y - Q$ are generated by local
 relations . The latter correspond to the singular points of the
 fibers of $p_h$.  If $f(Q_j) = P_i$ then the corresponding relation
 is described by $s_i^{N_i}$ as it was shown in section 2.1.  It
 finishes the proof of the theorem.

$\Box$

 The above theorem suggests that we can construct a big class of
 Burnside type groups as fiber groups.  We can consider the set of
 geometric vanishing cycles $s_i$ as a set of simple curves on one
 copy of the fiber.

The following construction allows to approximate  any 
topological data by the geometric ones.

Let $M_{g}^{L}$ be the compactified moduli space of curves of genus
$g$ corresponding to a subgroup of finite index $M_L$ in the group
 $Map(g)$. It is an algebraic variety
with quotient singularities only 
which contains a family of 
 similar type irreducible  divisors $S_I$ with normal crossings 
corresponding to different type of stable degenerations.
Singular points of $M_{g}^L$ correspond to stable curves with
 automorphisms and constitute a subset of codimension more than $1$
if genus $g > 2$. Let $M_g^{L_0}$ be an open nonsingular subvariety
in $M_g^L$ which corresponds to smooth curves of genus $g$ without
extra automorphisms. Then $\pi_1(M_g^{L_0}) = M_L$ .

The space $M_g^0 , g > 3$ is far from being affine. If fact there
is a natural 
map of $M_g^L$ into the  Satake compactification of the
moduli of  abelian
varieties of dimension  $g$ with a principal polarization. 
It maps each stable curve to
a point   corresponding to the Jacobian
 of its normalization. Thus all  divisors
corresponding to degenerate curve  have 
images of codimension at least $2$ if $g > 2$ under this map.
In particular generic hyperplane sections of Satake compactification produce complete  curves which lie in $M_g^0$.

We are interested in constructing  holomorphic families of
curves with a given set
of singularities.
 The following construction shows that there almost no restrictions in constructing such families over a disc.

\begin{lemma} Let $X_0$ be curve of genus $g > 2$ with a given
set of smooth  noncontractible cycles  $s_i^k$ on it.
 Assume that for a given $k$ all cycles $s_i^k, s_j^k$ don't intersect and correspond to different conjugated classes in the
 fundamental group of $X_0$. Then there is a holomorphic family
of curves over a disc $D$ which contains $X_0$ as nonsingular
fiber, the cycles  $s_i^k$ correspond to the vanishing
 cycles for degenerate fibers and monodromy group is generated by
 the products of Dehn twists over $s_i^k$ for each $k$.
\end{lemma} 

{\bf Proof} The condition on $s_i^k$ means that each set 
$s_i^k$ corresponds to some type of stable degeneration $I(k)$
modulo the action of $Map(g)$. Consider  a small
 complex disc $D_k$ around a generic point of $S_{I(k)}$ in $M_g$.
There is a local family of stable curves over $D_k$ which consists
of smooth curves outside the point of intersection of $D_k$ and
$S_{I(k)}$.
Let $p_k$ be a point on the boundary circle $dD_k$.
We can find a path $t_0^k$ connecting $0$ and $D_k$ inside $M_g^0$
which provides with a diffeomorpisms $X_0 \to X_{p_k}$ which
 maps $s_i^k$ into a family of vanishing cycles on $X_{p_k}$.
Indeed different paths provide with maps which differ by
the elements of $Map(g) = \pi_1(M_g^0)$.

Let us take an extension of $t_0^k$ into  smooth real analytic
 curves which ends up transversally 
at the intersection point of $D_k$ and $S_{I(k)}$. We can assume that  all the curves $t_k^0$ meet at $0$ 
being tangent to some one-dimensional complex subspace.
Thus we constructed a one dimensional ``octopus'' $W$ consisting
 of the extended curves $t_k^0$. After a small variation we can
complexify the resulting one-dimensional real set into a complex
disc $D$ which contains $0$ and intersects a given set of divisors
$S_{I(k)}$ with a prescribed monodromy corresponding to the
 connecting path $t_k^0$. It follows from the fact that a small
neighborhood of $W$ is Stein and 
holomorphic functions on it approximate continuous functions on $W$. The family of stable curves induced on $D$ satisfies the lemma.

$\Box$

\begin{lemma} For any finitely presented group $\Gamma$ we can construct a relatively projective family of curves $X_t$  over a holomorphic disc which has $\Gamma$  as a fundamental group.
\end{lemma}

{\bf Proof} For any group $\Gamma$ above we can find a surjective
homomorphism $r : \pi_g \to \Gamma$ for some $g$. Let
$N$ be the kernel of $r$. It has a finite number of generators $k_i$
 as  normal subgroup of $\pi_g$. The  elements $k_i$ can be realized
 as  cycles with normal intersections
 (including  selfintersections ) only on the
 curve $X_g$. Let us add a handle at each intersection.
 We can lift $k_i$ in a new Riemann surface $X_h$ into  a family
 of cycles $\tilde k_i$ without selfintersections.
 Let us  add to this family of cycles the generating cycles $a_j,b_j$
 of the additional  small handles. We obtain the set of conjugated
classes
 which
generates the kernel of the projection $\pi_h \to \Gamma$ and
each of the elements in this set
 is realized by a smooth cycle.

 Introduce a complex structure on $X_h$ and  consider positive 
 Dehn twists corresponding to $\tilde k_i, a_j, b_j$.
   All these cycles correspond to the conjugation classes  which
 belong to the kernel of the projection 
 $r_h : \pi_h \to \Gamma$. Note that
 the Dehn
 twist along the
 cycle  $s$  acts trivially on the quotient of $\pi_h$ by a normal
 subgroup generated by $s$. Therefore we can apply lemma 3.1 to
 $ X_h , \tilde k_i,a_j,b_j$ and obtain a relatively projective family $X$ of curves of genus $h$ over a disc  with $\pi_1 (X) =\Gamma $.

 $\Box$ 

Though the holomorphic families of curves look very 
 different from algebraic ones  we can rather easily
 embed a small deformation of such a family  into an 
 algebraic one.

Denote by $M$ the monodromy group of the family of curves over disc
 described in the lemma.
Note   that $D$ is a Stein subvariety in $M_g$. We can lift
it into any covering
 of $M_g$ which is unramified along $D$. 
  In particular we can lift $D$ into
 any variety $M_g^L$ for a subgroup of finite index
$M^L\subset Map(g)$ containing $M$.
 Since $D$ lies in affine subset of $M_g$ we can find
 an algebraic
 curve $C\subset M_g$( respectively
 $ M_g^L$) which  contains 
  a small variation of $D$.
 The family over $C$ induced from $M_g$ or
$M_g^L)$ extends
 globally a small variation of the family over $D$ without
 changing its topological data.
By taking generic $C$ we can assume
that the resulting family  $Y$ of curves  over a normalization of
$C$ is 
 smooth and the projection $p: Y\to C$ is a Morse type map. Thus having a family over disc with arbitrary data 
$ g, Ms_i$ where $M$ is generated by the Dehn twists over $s_i$ we
 obtain the 
 following  geometric data $ (g, M_Ls_i, M_L s_j)$ for any
 subgroup of finite index in $M_g$ containing $M$.
Accordingly considering any data $ g,Ms_i $ we can obtain 
a Burnside type approximation $ (g, M_L s_i^{N_1}, M_L s_j^{N_2})$
for any $N_1, N_2$.
By taking $N_1 = N$ and $N_2 = N^B$,  increasing the integer  $B$
and decreasing $M_L$ we obtain
series which approximate the group  given by $(g, M s_i^N)$.

\begin{rem} It seems plausible that most of the
 groups in such a series are 
 nonresidually finite and 
 violate any other good properties of  linear groups if  
 the group $ (g, Ms_i)$  violates them.

\end{rem}

\begin{rem} The above construction can be applied also to 
algebraic manifolds parametrizing special  curves instead of  moduli spaces. It results in different series of
groups as monodromy groups.
\end{rem}

\subsection {Potential examples}

  Now we suggest a construction that can lead to rather simple
 new examples of surfaces with nonresidually finite groups.  We thank
 V. Alexeev ,  S. Keel and M.Nori , for the fruitful discussions of the construction.

 Consider the map $f$ of the moduli space $M_g^L$ into 
 Satake compactification $S_g$ of the moduli space of abelian 
 varieties of dimension $g$ with a principal polarization.
 The latter is a projective variety which has a representation as a union
of $A_g, A_{g-1}..... $ where $A_g$ is a quotient of 
 the  space of positively defined hermitian matrices of
 rank $g$ by the action of $Sp( 2g,\Bbb{Z})$.
 If $x\in M_g^L$ corresponds to a stable curve $X$ then $f(x)\in S_g$ is
 a point corresponding to the Jacobian of the normalization of the
 curve $X$. 
Denote  by $SM_g$ the closure of the image $f(M_g)$ in $S_g$.
If $g > 3$  the map $f$ contracts analytic subvarieties corresponding
to the degenerate curves and curves with nontrivial automorphisms
into proper analytic subvarieties in $SM_g$ of codimension at least
 two.
 Denote by $\Delta_0$ the divisor in $M_g^L$ corresponding 
 to irreducible stable curves with one node.
 The image of $f(\Delta_0)$ consists of all the points of
$S_{g-1}\subset S_g$ which correspond to the Jacobian varieties
of dimension  $g-1$. 
Generic jacobian is known to be a simple abelian variety.
On the other hand the images of divisors corresponding to
 different
 type of stable degeneration intersect $S_{g-1}$ in proper 
subvarieties corresponding to nonsimple abelian subvarieties
(they decompose into a product after isogeny).

Consider the  surface $V$ in $SM_g$ obtained by a set
 of hyperplane sections.

 We can assume that :

1) The fundamental group of an open part of $V$ surjects
onto $Map(g)$.

2) The surface $V$ intersects the images of different divisors
$S_{I}$ at a
point only.

3) The intersection of $V$ and $f(\Delta_0)$ does not include
points from the images of other divisors $S_I$ or 
 subsets corresponding to curves with automorphisms.

 Denote the latter as $R \subset V$.
 By resolving $V$ at $R$ only  we again obtain a projective 
 surface $V'$. General hyperplane section $C$ of $V'$ will
 lift into a curve $C'$ in $M_g^L$ which intersects
only $\Delta_0$ and $C' - \Delta_0$ is contained in
$M_g^{L,0}$.
Thus we have 
  a family of curves $X_g$ over $C'$
 which has singular fibers of one type only and the
 monodromy group of it coincides with $M^L$.

Now we apply our construction to the above surface $V'$ and get a surface $V^{N}$ whose fundamental group $\pi_{1}(V^{N})$ is a  of a finite index in an extension of the fundamental group of a Riemann  surface by $\pi^{g}/(x^N = 1)$. Therefore  Zelmanov's
 conjecture implies that  the group $\pi_{1}(V^{N})$ is nonresidually finite. Hence we obtain a series of
 simple potential examples of  surfaces with 
nonresidually finite
 fundamental groups.

 The considerations from the previous subsection allow us to get even bigger variety of examples. Let us make a:

\begin{defi} Let $x$ be a primitive element in the fundamental group $\pi_{1}(g)$ of a Riemann surface of genus $g>1$ and $M^{L}$ be a subgroup of finite index in $Map(g)$. Consider the orbit of $x$ under 
$M^{L}$, $(M^{L} x)$ and take the $N$-th powers of all this elements. Consider the normal closure of this powers in  $\pi_{1}(g)$. Let us denote this normal closure by $(M^{L} x)^{N}=1$. We will denote the quotient group by $\pi_{1}(g)/(M^{L} x)^{N}=1$. (
Observe that definition depends on the choice of $x$.)

\end{defi}

Now the following generalization of the conjecture of Zelmanov's gives us a way of constructing more examples of surfaces with nonresidually finite fundamental groups.

{\bf Question} (Zelmanov)   For big $g$ and $N$ the groups  $\pi_{1}(g)/(M^{L} x)^{N}=1$ are nonresidually finite for any primitive $x$ and  $M^{L}$  a subgroup of finite index in $Map(g)$.

\section{ Some remarks on Shafarevich's conjecture for fibered surfaces}

\subsection{ The case of projective surfaces}

In this section we consider potential counterexamples to the
Shafarevich's conjecture based on our construction.  We begin with the general setting.  Let $f:X\to R$ be a Morse type
fibration with $X_t$ as generic fiber. Suppose that the fiber $X_0$ is
singular and has more than one component $X_0 = \cup C_i$. We also
assume that all components $C_i$ are smooth and without
selfintersection.  Denote the intersection graph of $X_0$ by
$\Gamma_0$.  Consider the retraction $cr : X_t\to X_0$ of generic
fiber on special fiber( see section 1).

 The preimage $cr^{-1}(C_i) $ in $X_t$ is an open Riemann surface with
 a boundary consisting of geometric vanishing cycles corresponding to
 the intersection points of $C_i$ with other components of $X_0$.  The
 fundamental group $\pi_1(cr^{-1}(C_i)$ is free. The natural embedding
 $cr^{-1}(C_i)$ into generic fiber $X_g$ defines an embedding of the
 fundamental groups $\pi_1(cr^{-1}(C_i)) = \Bbb{F}_i\to
 \pi_g$. Similar construction holds for any proper subgraph of curves
 in $X_0$.

\begin{defi} For any proper  subgraph
 $ K\subset \Gamma_0$ define a subgroup $ \Bbb{F}_K \subset \pi_g$ as
a fundamental group of the preimage $cr^{-1}(\cup C_i) ,i\in K$.
\end{defi}

\begin{rem} If the graph $K$ is connected then its preimage in
$X_t$ has only one component and vise versa.
\end{rem}

 Let $\pi_{1,f}$ be a fiber group obtained from $\pi_g$ by our base
 change construction for some $N$.

\begin{lemma} Suppose that there is a decomposition
 of a connected subgraph $K\subset \Gamma_0$ into a union $K_1 \cup
 K_2$ so that the image of $\Bbb{F}_K$ in $\pi_{1,f}$ is infinite, but
 the image of both $\Bbb{F}_{K_1} , \Bbb{F}_{K_2}$ is finite then the
 Shafarevich conjecture is not true.
\end{lemma}

{\bf Proof} Indeed under the conditions of the lemma we obtain an
infinite connected graph of compact curves in the universal covering
of the surface $S^{N}$.

$\Box$

Now we choose $N = 3$. In this case we can apply the above lemma due
to the  group theoretic result which concerns the quotients
of the free groups.

\begin{defi} Let $\Bbb{F}^g_k$ be a free group with $k$ generators
with a realization as a fundamental group of curve minus one or two
 points ( depending on $k$). Define $P^g( k,3)$ as the quotient of $\Bbb{F}^g_k$ by the
 set of relations $x^3 = 1$ for all primitive elements of 
$\Bbb{F}_k$ which can be realized by smooth nonintersecting curves in the
 above geometric realization 
 $ \Bbb{F}^g_k$ of $\Bbb{F}_k$.
\end{defi}

\begin{theo} The group $P^g (k  ,3)$
 is equal to the Burnside group $B(k ,3)$ and hence finite.
\end{theo}

{\bf Proof} See Appendix B.  $\Box$

 \begin{corr} The group $\pi_g/(x^3 = 1)$ is a quotient of $B(2g,3)$
 by one additional relation.

\end{corr}

 Let $X_0$ be a graph of smooth curves $C_i$ with each curve
intersecting at most two others ( chain or ring). 
 Suppose that there
exists a Morse family of curves $X \to R$ with a set of vanishing 
 cycles $VS$ and  a monodromy group $M$ which has $X_0$ as
 fiber.
Assume that  cycles from $VS$ correspond to different
 singular points of $R$ unless they correspond to $X_0$.
Assume that  $VS$ is  decomposed into a union $S_0\bigcup S_1\bigcup S_2$. 
 Assume that  cycles from different subsets $S_i,i=0,1,2$ correspond to different
 singular points of $R$.

\begin{lemma} Suppose that the monodromy group $M$ and the sets
of vanishing cycles $S_0,S_1,S_2$ satisfy the following properties:

 1.  The image of $\pi_1(cr^{-1}(C_i)) $ in $\pi_g/ Ms_j^3  $ is finite
 for any $i, s_j \in S_0\bigcup S_1$.

 2. The quotient group $\pi_g / ((Ms_j^3), Ms_k^{3^S}) $ is infinite for some integer  $S$ and $s_j \in S_0\bigcup S_1 , s_k\in  S_2 $.
Then the universal covering $\tilde X$ is holomorphically nonconvex.
\end{lemma}

{\bf Proof} Indeed the universal covering $\tilde X$ contains 
an infinite covering of $X_0$ which is connected and consists 
of compact curves. Hence $\tilde X$ is not holomorphically convex.

$\Box$ 
 
We are going to construct a family of surfaces which presumably contain an infinite number of surfaces with above property. We would like to produce such families from a standard 
family over an interval $I$.

Let $g > 3$ and consider a family over an interval $I = [0,1]$ which
has a fiber $X_0$ over $0$ and a singular fiber $X_1$, 
The generic fiber $X_t$ surjects on $X_0$ and $X_1$ and 
vanishing cycles for both singular fibers are realized as
 smooth curves on $X_t$. Denote the corresponding set of
cycles as $S_0, S_1$ respectively.
 We assume that they don't intersect
and for any component  $C_i\subset X_0$ there is a
 corresponding cycle $s_i\in S_1$ which projects into a smooth
homologically nontrivial cycle in $X_i$.

This family over interval can be complexified 
 into an algebraic family over
a complete curve $R$.
We can assume that the monodromy of the resulting family
is any subgroup of finite index in the group $Map(g)$
 which contains commuting monodromy transformations $T_0,T_1$
 defined by the fibers $X_0, X_1$.

Let $M_D$ be a subgroup of $Map(g)$ which commutes a
 contraction map $\pi_g\to \pi_1(X_0)$.
 Since the group $M_D$ contains $Map(g)$ for any curve $C_i$ we have obtain that the image $\pi_1(cr^{-1}(C_i))/(M_D s_i^3) $ is a
 finite group of exponent $3$ for any component $C_i$.

\begin{lemma} The quotient $\pi_g/ (M_D s_i^3)$ is infinite
if the number of components $C_i, g(C_i) > 0$ is more than $1$.
\end{lemma}

{\bf Proof} Indeed the group above maps surjectively
 onto a free product
of nontrivial  Burnside groups of exponent $3$ corresponding
to  different components of $X_0$ with  nonzero genus.
The latter is infinite which implies the lemma. 
$\Box$

Let $M^f$ be any subgroup of $M_D$ which contains $T_0,T_1$ and
 has  the  property that
the image of $\pi_1(cr^{-1}C_i) $ in $\pi_g /(M^f s_i^3)$ is finite.

The group $M^f$ defines a set of  
 subgroups $M_j^f$ of finite index in $Map(g)$ containing $M^f$.
The intersection of this set of subgroups  coincides with $M^f$ since 
$Map(g)$ is residually finite.
  
For each $M_j^f$ we can find a curve $R_j$ with a  Morse 
family of curves
$X_j$  which contains a
topological family over interval constructed above and with
 a monodromy $M_j^f$. Let $S_2^j$ be a complementary set of vanishing cycles in the family $X_j$. We can now consider any 
$S$ and construct a new family $X_j^S$ using the theorem 3.1
 with $N_0 =N_1 = 3 ,N_2 = 3^A$.

\begin{corr} For  $M_j^f $ and integer $A > 1$ we obtain a group
$\pi_g /((M_j^f)s_i^3 , M_j^f s_k^{3^A})$ as a fiber group.
Here $s_i \in S_0\subset S_1 , s_k\subset S_2$.
 
\end{corr}

\begin{rem} The subset $S_2$ depends on the actual curve
$R_j$. The dependence of the
 fiber group on $S_2$ weakens with $A$ converging  to infinity.
The resulting family of groups  approximates the infinite
 group  $\pi_g/ (M^f s_i^3)$ as $M_j^f$ converges to $M^f$ and
 $A$ converges to infinity.
\end{rem}

\begin{con} Let $X_0$ be a nontrivial chain or ring  of curves
of genus greater than $0$. Assume that $g$ is an arithmetic
genus of $X_0$ and $X_I$ is  a family of curves 
over an interval described above. 
 For any small enough  subgroup $M^f\subset M_D$
 defined above there exists a subgroup $M^f_j\subset Map(g)_f$ 
of finite index and  an integer
$A > 1$ such that $\pi_g/ (M^f_j s_i^3, M_j^f s_k^{3^A})$ is infinite
for $s_i\subset S_0\bigcup S_1$ and any finite subset $s_k\in \pi_g$
\end{con}

 If the answer to  the above conjecture is positive 
 then the Shafarevich
 conjecture is not true.
 The fact that the family of groups  parameterized  by $M_j^f$ and $A$
approximates a group which has a nontrivial  free product of groups
as a quotient provides with a strong evidence supporting the above
conjecture. On the other and the  resulting fiber group does not have infinite
linear representations which are
 equivariant with respect to the action of the monodromy group  ( \cite{LN}).

There is also another possibility to satisfy condition 2 .  We can
easily construct a family of curves of genus $g=2k, k>1$ such that
there is a fiber in this family which consists of two components each
of genus $k$ that give us a tree of components.  Applying the base
change construction for a given $N$ we obtain a surface $S - Q$ with
the image of the fundamental group of every component in $\pi_1(S -
Q)$ being equal to $P^g(2k,N) $. The fiber group of $S-Q$ is equal to
$\pi_{2k }/(x^{N}=1)$. Now as we have shown we have the same behavior
on $S^{N}$ for the closed curves and surfaces. We can formulate the
following question:

{\bf Question} Are there such a $k$ and $N$ such that $P^g(2k,N)$ is a
 finite group and $\pi_{ 2k }/(x^{N}=1)$ is an infinite group?

If the answer of the above question is affirmative for some $ N, k$ we
  get a counterexample to the Shafarevich conjecture. We should point
  out that if the groups obtained from the components are finite then
  $\pi_{2g}/(x^{N}=1)$ does not have infinite linear representation. If the
  Shafarevich conjecture is correct the answer of the above question
  is negative. It also implies the answer to many similar group
  theoretic questions. The most basic question seems to be the
  following:

{\bf Question} Is there such an $N$ and such $2\le m_{1}< m_{2}$ for
 which $B(m_{1},N)$ is finite and $B(m_{2},N)$ is infinite?

Recently we were informed by Zelmanov that he can show that there
 exists an integer $d(0)$ so that for every prime number $p$ and an
 integer $d$ the group $B(d_{0},p)$ is finite if and only if the group
 $B(d,p)$ is finite. This result suggests the existence of an abstract
 group theoretic version of holomorphic convexity.  The above
 considerations indicate a possibility for analysis of infinite groups
 by analytic methods.

\subsection{Other applications}

In this subsection we discuss symplectic and arithmetic versions of
 our construction.  As Gompf has shown \cite{GOMPF}, every finitely
 presented group can be realized as a fundamental group of a
 symplectic manifold.  It is reasonable to ask the if being symplectic
 puts any restrictions on the structure of universal coverings. Our
 construction easily extends to the symplectic category. It leads to
 interesting examples of symplectic four dimensional manifolds ( see
 \cite{BKS}).  Using this construction we have defined in \cite{BKS}
 an obstruction to a symplectic Lefschetz pencil being a K\"{a}hler
 Lefschetz pencil.

We begin with a symplectic fourfold $X$.  Consider the corresponding
 Lefschetz pencil with reducible fiber and apply to it the our
 construction. So we get for a fix integer $N$ a symplectic fourfold
 $S^{N}$.  Let $\rho$ be a generic representation $\rho :\pi_{1}(S^{N}) \to
 GL(n,\Bbb{C})$ whose image is not virtually equal to $\Bbb{Z}$.  Denote by $Y_{i}$ the components of the preimage of
 the reducible fiber of $S$ in $S^{N}$ and denote by $F$ the general
 fiber of $S^{N}$. Denote by $\Gamma$ the image of $\pi_{1}(F)$ in
 $\pi_{1}(S^{N})$ and by $\Gamma_{i}$ the images of the fundamental
 groups of $Y_{i}$ in $\pi_{1}(S^{N})$.

If the restrictions of $\rho$ on $\Gamma$ and $\Gamma_{i}$ for all $i$
 are both finite or infinite we will say that the obstruction
 $O(X)^{N,n}$ is equal to zero and to one otherwise.

\begin{prop}

If $S$ is K\"{a}hler $O(X)^{N,n}$ is trivial for every pair $N,n$.
\end{prop}

Indeed otherwise we contradict the holomorphic convexity of the
covering corresponding to infinite linear representation of the
fundamental group. In \cite{BKS} we have constructed examples of
symplectic fourfolds with nontrivial $O(X)^{N,n}$.

In this article we consider only surfaces with given projection on a
curve. Algebraically this means that the field of rational functions
on the surface is provided with a structure of a one-dimensional field
over a field of rational functions of the base curve. This picture is
parallel to that of a curve defined over number field.  It suggests
that there should be a natural arithmetic version of the Shafarevich's
conjecture. The notion of holomorphic convexity is not well defined in
the arithmetic case. Instead we can describe the analogue of the
absence of infinite chains of compact curves in the universal covering
in the arithmetic case. The absence of infinite chains of compact
curves seems to be the only obstruction to the holomorphic convexity
in the case of compact complex surfaces.

 Let $C$ be a projective semi-stable curve over $K$.  Since $K$ is not
algebraically closed we obtain a nontrivial map $C\to Spec (O_K)$
where $O(_K)$ is the ring of integers.  Extending $K$ if necessary we
can assume that $C$ is a semistable curve.  That means $C$ is a normal
variety with semi-stable fibers consisting of normally intersecting
divisors.  The function field $K(C)$ is regular and has dimension one
over $K$.  Any maximal ideal $\nu$ in $O_K$ defines a subring
$O_{\nu}$ of $K(C)$ which consists of the elements of $K(C)$ which are
regular at the generic point of any component of the preimage of
$\rho$ in the scheme of $C$. The ring $O_{\nu}$ contains the ideal of
elements which are trivial on the preimage of $\nu$.  Denote this
ideal by $M_{\nu}$. The quotient ring $O_{\nu} / M_{\nu}$ is a finite
sum of fields of rational functions on the components of the fiber of
$C$ over $\nu$.  Consider the maximal nonramified extension
$K(C)^{nr}$ of $K(C)$.  It is a Galois extension with a profinite
Galois group $Gal^{nr}(K(C)$.  The field $K(C)^{nr}$ contains the
maximal nonramified extension $K^{nr}$ of the field $K$ as a
subfield. (Observe that the field $\Bbb{Q}^{nr} = \Bbb{Q}$ but for
many other fields $K$ the field $K^{nr}$ is an infinite extension.)
The group $Gal^{nr}(K(C)$ maps surjectively onto $Gal^{nr}(K)$.
Denote the kernel of the corresponding projection as $Gal_g^{nr}K(C)$.

 Any maximal ideal $\rho$ in the ring of integers $O_{K^{nr}}$
contains the unique maximal ideal $\nu$ of $O_K$. We can now define
the subring $A_{\rho}$ as an integral algebraic closure of the subring
$O_{\nu}$ in $K^{nr}(C)$. Let $Res_{\rho} = A_{\rho}/I(\rho)$ be the
quotient ring by the ideal generated by $\rho$ in $A_{\rho}$. It is a
semisimple ring of finite characteristics.  We formulate a strong
arithmetic analogue of the Shafarevich's conjecture.

\begin{con} Let $K(C)$ be a field as above. For some finite extension
 $F$ of $K$ we can find a semistable model $C'$ of the field $F(C)$
such that the ring $Res_{\rho}$ is a direct sum of a finite number of
fields for any ideal $\rho$ of the field $F^{nr}$.

\end{con}

 We can also formulate this conjecture in a more geometric
 language. Namely if $C'$ is a semistable curve over $F$ then for any
 finite nonramified extension $L$ of $F(C')$ we have a model $C_L$
 with a finite map onto $C'$. The fibers of this new model $C_L$ are
 uniquely defined by $C'$.

\begin{con} There is a number $J(F(C'))$
such that the number of components of any fiber of $C_L$ is bounded by
$J(F(C'))$ for any finite extension $L$ of $F(C')$ containing in
$F(C')^{nr}$.
\end{con}

\begin{rem}The result of the conjecture depends ( at least formally)
 on the chosen model $C'$.  If we blow up a generic point on a fiber
the conjecture becomes false if the corresponding covering induces an
infinite covering of the fiber.  Thus the arithmetic version is
sensitive to the change of the semistable model whereas the geometric
conjecture is not.

\end{rem}

\begin{rem} The fields $F(C')^{nr}$ correspond to the factor of $\pi_{1}^{fin}(C')$ that is acted trivially by all inertia groups.
 ( Here we denote by $\pi_{1}^{fin}(C')$ the profinite completion of the geometric fundamental group of $C'$ ). 
 Geometrically this means that we consider coverings of $C_L$ that are unramified over a generic point of every irreducible divisor in $C_L$.

\end{rem}

There exists some evidence for the above arithmetic conjectures.
 Partial results in this direction were obtained by the first author
 several years ago (1982).  Namely he proved that the torsion group of
 an abelian variety $A$ is finite for any infinite algebraic extension
 of $K$ which contains only finite abelian extensions of $K$.

 In particular it is true for
 infinite nonramified extension of $K$ where $K$ is a finite extension
 of $\Bbb{Q}$. (This
  result was announced at Delange-Puiso seminar in Paris, May
 1982 and later appeared in \cite{Co}).
  Yu. Zarhin  proved that the same is true if the
 infinite  extension of $K$ contains only
 finite number of roots of unity under
 some conditions on the algebra of endomorphisms $A$.

The above result states that the group $Gal_g^{nr}K(C)$ has a finite
 abelian quotient and hence the conjecture 4.2 is evidently true for
 the quotient of  $Gal^{nr}K(C)$ by the commutant of
 $Gal_g^{nr}K(C)$. For the same reason it is true for the
 quotient of $Gal_g^{nr}K(C)$ by  any iterated  
commutant.

\begin{rem}
  Any curve over arithmetic field can be obtained as a covering of
$\Bbb{P}^1$ ramified at $ (0,1,\infty) $ according to the famous
Bely's theorem. As it was pointed out by Yu.Manin any arithmetic curve
has a ramified covering which is a nonramified covering of a modular
curve. Therefore one might reduce the above conjecture to modular
curves being considered over different number fields.  In the complex
case we don't have similar simple class of dominant manifolds( see the
discussion in \cite{BH}.
\end{rem}

If Conjecture 4.2 is correct for every curve over every finite
extension of $\Bbb{Q}$ then we get that there isn't any infinite
chains of compact curves on the universal covering of a projective
surface with a residually finite fundamental group. According to a
conjecture of M. Ramachandran this is the only obstruction to
holomorphic convexity for the universal coverings of a projective
surfaces.

\section{ Appendix A - Fiber groups and monodromy}

This appendix includes several technical results from the theory of
 surfaces which we use in our article. First theorem generalizes the
 result we used in order to pass from quasiprojective surface to the
 projective smooth surface in section 2.

\begin{theo} Let $V$ be a normal projective
 surface with $Q$ being the set of singular points in $V$ and $f: V'
 \to V$ be a finite surjective map from another normal projective
 surface $V'$.  Assume that for any $q\in f^{-1}(Q)$ the map of the
 local fundamental groups $ f_* :\pi_1( U_q - q) \to \pi_1(V - Q)$ is
 zero, where $U_q$ is small topological neighborhood of $q$.  Then
 there is a natural map $ f_* : \pi_1(V') \to \pi_1(V - Q)$ which is a
 surjection on a subgroup of finite index bounded from above by the
 degree of $f$.

\end{theo}

{\bf Proof} Indeed $Q$ is a set of isolated singular points.
 Following the proof of the lemma 2.5 we obtain a map $f_* : \pi_1 (V'
 - f^{-1}(Q)) \to \pi_1(V - Q)$. In order to prove the theorem it is
 sufficient to show that the kernel of natural surjection $i_* :
 \pi_1(V' - f^{-1} (Q))\to \pi_1(V')$ lies in the kernel of $f_* $.
 The preimage $f^{-1}(Q)$ also consists of a finite number of points
 and the kernel of $i_*$ is generated as a normal subgroup in
 $\pi_1(V' - f^{-1}(Q)$ by the local subgroups $\pi_1(U_q - q), q \in
 f^{-1}(Q)$. Due to the condition of the theorem the images of these
 groups are trivial in $\pi_1(V - Q)$ and hence we obtain the map
 $f_*$.  $\Box$

We shall also need the following general result which allows to
 compare the properties of fiber groups and fundamental groups of the
 smooth quasiprojective surfaces.

 \begin{theo} Let $ f : X \to R $ be a quasiprojective surface which
 has a structure of a family of curves over a smooth curve $R$. Assume
 that generic fiber is an irreducible smooth projective curve and any
 fiber contains a component of multiplicity one.  Let $h: \pi_{1,f}
 (X) \to G$ be $M_X$ invariant homomorphism into a finite group $G$
 with $ M_X$ action.

 Let $K$ be the kernel of $h$. Then there exists a subgroup of finite
index $ H \in \pi_1(X)$ that the intersection of $H$ with
$\pi_{1,f}(X)$ is equal to $K $.
\end{theo}

\begin{rem} This is true if there exists a
section $s:\pi_1(R) \to \pi_1(X)$, since we can define $H$ as a
  subgroup of $\pi_1(X)$ generated by products of the elements of $
  s(\pi_1(R))$ and $ K$. If the curve $R$ is open then the group
  $\pi_1(R)$ is free and hence a section always exists.

 \end{rem}

{\bf Proof} The group $K$ is a normal subgroup of $\pi_{1}(X)$ since
 it is invariant under the conjugations from both $\pi_{1}(R)$ and
 $\pi_{1,f}(X)$.  Denote by $Q$ the quotient $\pi_{1}(X) / K$. It is
 an extension of $\pi_{1}(R)$ by $G$. Hence there is an action of
 $\pi_{1}(R)$ over $G$. Since $G$ is a finite group $\pi_{1}(R)$
 contains a subgroup of finite index which acts trivially on $G$ via
 interior endomorphisms.  This subgroup corresponds to a finite
 nonramified covering $ \phi : C \to R$ and will be denoted by
 $\pi_1(C)$. Its preimage $Q'$ in $Q$ is a subgroup of finite
 index. Consider a subgroup $K_G$ of $Q'$ consisting of the elements
 commuting with all the elements of $G$.  The group $K_G$ projects
 onto $\pi_1(C)$ by the definition of $\pi_1(C)$ and therefore its
 intersection with $G$ is a cyclic subgroup of the center of $G$.  The
 group $\pi_1(C)$ contains a subgroup of finite index $K'$ where the
 corresponding central extension splits.  Namely if the order of the
 corresponding cyclic extension is $n$ then it splits over any cyclic
 covering of order $n$.  The preimage of $K'$ in $Q'$ splits into a
 direct product of $G$ and $K'$. Hence on the preimage of a subgroup $
 K'\subset \pi_1(C)\subset \pi_1(R)$ we have a natural extension of
 the projection from $\pi_{1,f}(X) \to G$. $\Box$

\section{ Appendix B - Some group theoretic results}

This appendix contains several  group theoretic results.
We present proofs though  
most of the results are presumably known to the experts
in the group theory..

Recall that we denote by $BT(n,m)$ the quotient of a free group
 $\Bbb{F}^n$ by a normal subgroup generated by the elements
 $x^m = 1$ where $x$ runs through all primitive elements of
 $\Bbb{F}^n$.

\begin{prop} If $m$ is divisible by 4 and $n \ge 2$ then $BT(2,m)$
is infinite.
\end{prop}

{\bf Proof}  The group $BT(2,4)$ has an infinite
 representation. Namely let $\Bbb{Q}_8 $ be the group of the unit
 quaternions of order $8$. It acts on $ H =\Bbb{R}^4$ by
 multiplications.  Consider a group of the affine transformations of
 $\Bbb{R}^4$ generated by two generating rotations $g_1, g_2$ in
 $\Bbb{Q}_8$ but with different invariant points.  The resulting group
 $G$ will be infinite. It has two generators and $G$ has $\Bbb{Q}_8$
 as its quotient group. Any element of $g \in G$ which projects
 nontrivially into $G$ has order $4$. Indeed $g$ is an affine
 transformation of $\Bbb{R}^4$ with a linear part of order $4$ and
 without $1$ as eigenvalue. Hence $g$ in the conjugacy class of  its linear part and has order $4$. 
 
$\Box$

\begin{rem} Similar argument can be applied to any finite
 subgroup of a skewfield instead of $\Bbb{Q}_8$ (see the description
 of such finite groups in J.Amitsur, Ann. of Math., 1955, vol. 62,
 p. 8). The above result makes it plausible that 
 the groups
 $P(n,m)$ are infinite if $ m > 3, n > 1$. 
\end{rem}

The case $BT(n,3)$ is different. The following result gives a hint on
the effects which occur with the exponent $3$.

\begin{lemma} Let $G_i, i=1,2$ be the groups generated by $ a,b$ with
 relations

1)  $G_1 : {a^3 = b^3 = (ab)^3 = 1}$

2) $G_2 : {a^3 = b^3 = (ab)^3 = (ab^2)^3 = 1}$.

Then $G_1$ is infinite and  $G_2$ coincides with the
 Burnside group $B(2,3)$ and hence it is  finite.
\end{lemma} 
{\bf Proof} The group $G_1$ has a natural geometric
 realization. Namely let us take $\Bbb{P}^1$ minus three points $p_i , i = 1,2,3$.
The fundamental group $\pi_1(\Bbb{P}^1 - p_i) =\Bbb{F}_2$.
If we impose relations $x^3$ on the elements which can be realized by
smooth curves in $\Bbb{P}^1 - p_i$ then we obtain the set 1).
Let us take a $\Bbb{Z}_3$ character $\chi$ of $\Bbb{ F}_2$ which is nontrivial
on $a,b,ab$. We obtain a covering of $\Bbb{P}^1 $ ramified over
 three points. Imposing the above relations corresponds to the
 completion of the curve. Hence the subgroup of $G_1$ which
 is the kernel of $\chi$ coincides with the fundamental group
of the corresponding complete curve.  
  Since
 the above curve is a torus the group $G_1$  is an extension of $
 \Bbb{Z}_3 $ by a free abelian group $ \Bbb{Z}+ \Bbb{Z}$.  The
 description of 2) immediately follows since it is the quotient group
 of the group in 1).
 The element $ab^2$ generates $\Bbb{Z}+ \Bbb{Z} $
as a $\Bbb{Z}_3$-module. Hence $(ab^2)^3$ generates $3(\Bbb{Z}+ \Bbb{Z})$
and the resulting group is an extension of $\Bbb{Z}_3$ by the
group $\Bbb{Z}_3 + \Bbb{Z}_3$.

$\Box$

We shall use the following standard notations.
Denote by $(a,b)$ the commutator of the elements $a,b$ and
we put a sequence of brackets to denote an element obtained
by iteration of the procedure.

The group $B(n,3)$ has a rather simple description.
 It is a 
  metabelian group with a central series of length three.
 The elements $(a,x)$
 where $x\in [B(n,3),B(n,3)]$ are in the  center of $B(n,3)$.
 In fact  $((a,b),c)$ varies
under permutation according to standard $\Bbb{Z}_2$ character of the
 group $S_3$ for any $((a,b),c)$.
 In particular $((x,r),x) = 1$ for any $x$.

\begin{lemma} Let $G$ be a finitely generated group
 with a given
 set $S$ of generators. Assume that 
  $((a,b), f)$ is invariant under any even permutation of $a,b,f$
 for any $a,b\in S, f\in S\bigcup (S,S)$
 where
 the latter denotes the set of pairwise commutators of the
 elements from $S$. Then $G$ is a metabelian group
 and it has a central series of length 
 at most $3$.
\end{lemma}

{\bf Proof} The proof closely follows the  proof from \cite{MS}.
Indeed we can write $((a, b),(c, d))= (a,(c,d)), b)$.
  We  deduce next that  the above expression is invariant 
under even permutations. The latter implies that it is equal to 1
and hence $((a,b), c)$ commutes with any element from $S$ and hence
lies  in the center of $G$.
The quotient of $G$ by the center  $Z\in G$ is also generated
by $S$ with equality $((a ,b),c) = 1$  for any $a,b,c\in S$
which means that
$(a,b)$ is in the center of  $G/Z$ for any $a,b\in S$.
That means $G/Z$ is a central extension of abelian group which
 finishes the proof.

 $\Box$

\begin{corr} Under the above conditions the commutant
$[G,G]$ is additively generated by the
elements $(a,b),((a,b)c),a,b,c \in S $.
\end{corr}

\begin{lemma} Assume that $S$ in the above lemma consists of
$n$ elements ,  $[S] , [S,S]$ consists of elements of order $3$
and there is a surjective map $p: G\to B(n,3)$. Then $p$ is an isomorphism. 

\end{lemma}

{\bf Proof} 
 The commutant $[G,G]$ is additively generated 
as $G^{ab}$-module by the elements $(a,b)$ under the above conditions. Hence it  is of exponent $3$ and the number of elements in
$G$ is not greater than the number of elements in the commutant of $B(n,3)$.
The abelian quotient $G^{ab}$ is isomorphic to $B(n,3)^{ab}$.
Therefore the number of elements in $G$ is not greater than
in $B(n,3)$ and a surjective map $p : G\to B(n,3)$ is an isomorphism.

$\Box$

\begin{prop} The group $BT(n,3)$ is finite and coincides with
the Burnside group $B(n,3)$ if $BT(3,3) = B(3,3)$.
\end{prop}

{\bf Proof}
 Indeed this is true for $n= 2$ as it was shown above ( see e.g. \cite{MS}). The group $BT(n,3)$ has a natural surjective map onto
$B(n,3)$ and hence we have to check the condition on $(a,b),c)$.
The latter  is enough to check  for the group
 with three generators.

$\Box$

\begin{lemma} $BT(3,3) = B(3,3)$.
\end{lemma}
{\bf Proof}

The group $BT(n,3)$ is obtained as an extension of $B(3,2)$ by $c$
and since $cx$ is a generator for any $x \in B(3,2)$ we
 obtain that $ (cx)^3 = 1$.
 Therefore $x^{-1}c x , c$ commute for
any  $x \in B(3,2)$ and the kernel $K$  of the projection
 $p_c :BT(3,3) \to B(3,2)$ is an abelian group of exponent $3$.
The sum   $x^{-1}c x + c  + x c x^{-1} = 0$ and we can easily
 deduce that
$K$ has $4$ generators as a $\Bbb{Z}_3$-space.
Therefore  the group $BT(3,3)$ has the same number of elements
as $B(3,3)$ and since there is a natural surjection
$ BT(3,3)\to B(3,3)$ they are isomorphic.

$\Box$

 \begin{lemma} Let $S$ be a finite set of generators of the group $G$. Assume that any subset of $4$ elements in $S$ generates a
subgroup of exponent $3$. Then $G$ is of exponent $3$.
\end{lemma}

{\bf Proof}
 The assumption implies that any three-commutator of  $a,b,c\in S$ lies in the center $C$
 of the group $G$ and the group satisfies  lemmas 6.2, 6.3.

$\Box$

\begin{lemma} Assume that $G$ has four generators $a,b,c,d$ that
 $G$ contains a set of subgroups of exponent $3$ which includes
  groups generated by triples of generators  $ a,b,c, d$  and also
 the groups
$a,b,(cd)$,  $ c,d, ( a,b)$.
 Then $G$ is of exponent $3$.
\end{lemma}
 
{\bf Proof} It follows from the above results that $((a,b)c)$
 and similar combinations are invariant under cyclic permutations.
We also have $((a,b)(c,d)) = (((c,d)a)b) = (((d,c)b)a) = (((a,b)c)d)
= (((a,b)d)c)$. Thus the value of the above commutator
 does not depend on any even permutation of the symbols. On the other hand it transforms into opposite
if we permute $(c,d)(a,b)$. Hence  all the above elements are equal to 1. This implies the lemma.

$\Box$

We are interested  in the groups which occur as the
 quotients of the fundamental groups of curves.  Namely we will study the groups $P^{g}(n,m)$ defined in section  four.
 We are going to use a representation of the fundamental
 group of a Riemann surface  as a subgroup of index $2$ in the group
 generated by involutions which was successfully used by J.Birman, M.Nori, W.Thurston , B.Wainrieb and many others.  Let $X$ be a curve of genus $g$. It can be represented as 
 a double covering of $\Bbb{P}^1$ ramified over $2g+ 2$ points and
$\pi_g$ is a subgroup of index two in the group generated
by $2g + 2$ involutions $x_i$ with additional relation that
the product of all involutions is an involution again.
 Similarly the group $\pi_1(X_g - pt)$ is realized as subgroup
 of index two in the group generated by $2g + 1$ involutions. The fundamental group of a curve minus two points is realized as
a subgroup of index two in the group generated by $2g + 2$
involutions without any additional relations. 
 
The groups above are free groups, but they are provided with a special realization as the fundamental groups of open curves.

\begin{defi} We shall  denote by $\Bbb{F}^{g}_{n} $ a free group of $n$ generators provided with a realization as a fundamental group of a curve of genus $g$ minus one or two  points. In case $n$ is odd then  $\Bbb{F}^{g}_{n} $ is realized as a fundamenta
l  group of a curve of genus $g$  minus one point. In case $n$ is even $\Bbb{F}^{g}_{n} $ is realized as a fundamental group of a curve of genus $g$  minus two points.
\end{defi}

Recall that $P^{g}(n,m)$ is a quotient of $\Bbb{F}^{g}_{n} $ by the relations $x^{m}=1$ where $x$ runs through the primitive elements of $\Bbb{F}^{g}_{n} $ 
which can be realized as a smooth loops in the above geometric realization.

The following lemma reduces a general case to $ n \leq 4$
 
 \begin{lemma} If $P^g(4,3) = B(4,3)$ then $P^g(n,3) = B(n,3)$
 for any $n\geq 4$.
\end{lemma}

 {\bf Proof} 

The group $\Bbb{F}_n^g$ is represented as subgroup of index two generated by involutions $ x_1,..,x_{n+1}$.
The set of standard generators of $\Bbb{F}_n^g$  can be taken as $ x_1x_i, i\neq 1$. Any four elements $x_1 x_j $ generate 
 a subgroup of $P^g(n,3)$ which is a quotient of $P^g(4,3)$
and by assumption of the lemma the latter is of exponent $3$. Hence by lemma 6.5 $P^g(n,3)$ is of exponent $3$.
 Since the set of generators includes only primitive elements $\Bbb{F}_n$
the group $P^g(n,3)$ coincides with $B(n,3)$.

 $\Box$

\begin{lemma} If $P^g(3,3) = B(3,3)$ then $P^g(4,3) = B(4,3)$.
\end{lemma}

{\bf Proof} The group $P^g(4,3)$ is obtained from the curve
 of genus 2 minus a point.
Consider a standard decomposition of $X^2$ into a union of two
 handles corresponding to pairs of generators $ a,b$ and $c,d$
of $\Bbb{F}^g_4$. There are topological embeddings of tori with two discs
 removed corresponding to the subgroups generated by
any three symbols from the set  $(a,b,c, d)$. The group $P^g(3,3)$ is realized as the 
quotient of the fundamental group of torus minus two discs. The above tori  are obtained from one of the
 handles by adjoining a neighborhood of a generator in another handle.  By assumption of the lemma 
 $((a,b),c)$ is transformed into itself or the opposite element
under the permutation of symbols.
Similar groups correspond to the triples $( a,b (cd))  ,( c,d, (ab))$.
Namely we can consider a corresponding handle minus a point. Thus we can apply  lemma 6.6 and obtain that $P^g(4,3)$
is of exponent $3$ if the group $P^g(3,3)$ is.

$\Box$

 \begin{theo} $P^g(3,3) = B(3,3)$.
\end{theo}

Let us first describe  the geometric picture.
  Consider the torus $T^2$ with a small
 embedded interval $I$.
 Let $p_1,p_2$ be two different points in the interval $I$
and  $I_1$ be the interval between $p_1,p_2$ inside $I$.
We assume that $ p_0$ is a point in  $I - I_1$ and identify $p_0$ as an initial point for 
 the fundamental group $\pi_1(T^2 - p_1- p_2)= \Bbb{F}^g_3$.
Assume that $T^2,I$ are given with orientation.
We consider smooth oriented loops through $p_0$ which are transversal to $I$.

\begin{lemma} Any element of $P^g(2,3)$ with $p_0$ as an
initial point  is represented by an oriented curve $A$
 without selfintersection which does not intersect $I$ and such
that
$A,I$ defines a standard orientation of $T$ at $p_0$.

\end{lemma} 
{\bf Proof} Let $a,b$ be standard  generators of $\Bbb{F}_2^g$ with
 a given orientation. The group $Out(\Bbb{F}_2^g) = SL(2,\Bbb{Z})$
 and can be realized by linear periodic map of torus. In particular we
 can represent topologically  the elements
of $SL(2,\Bbb{Z})$ by maps  which stabilize the points of $I$.
In this way we obtain any map of $\Bbb{F}_2^g$ into itself which
 transforms $aba^{-1}b^{-1} = C$ into itself.
Any homomorphism of the free group with above property
is induced by this action of $SL(2,\Bbb{Z})$.
Thus $g(a), g\in SL(2,\Bbb{Z}), a \in B(3,2)$ can be any element such
 that there exists $b'\in B(3,2)$ with $g(a)b'g(a)^{-1}b'^{-1}= C$
where $C$ is a given generator of the center.
But we can find such $b'$ for any $g(a)$ which is not in the center. The curve $g(a)= A$ will be the image of a map which is linear
 outside a neighborhood of $I$ and keeps orientation intact.

$\Box$

\begin{rem} Since $A^{-1} I$ represents the opposite orientation
the lemma actually shows that the element $A^{-1}$ can be represented by 
 another simple closed curve $B$ with orientation $B  I$ opposite to $A^{-1}  I$.
\end{rem}
 
The group $P^g(3,3)$ is generated by $B(2,3)$ realized as
 above and an element $r$ realized by a curve $R$ with one
 selfintersection. We have a natural representation
 of $r$ as $x_1 x_2^{-1}$ where
 $x_1, x_2$ are simple curves  with the same orientation which move
 around $p_1, p_2$ respectively. The element $c = x_1 x_2$ is a natural central loop in $P^g(2,3) = B(2,3)$.

\begin{lemma} Any element $bx_i $ is realized by a simple loop
if $b$ is not in the center of $B(2,3)$.
\end{lemma}

{\bf Proof} We have to find a simple representative of $b$
through $p_0$ with an appropriate
orientation. The latter exists due to the previous lemma.

$\Box$

\begin{lemma} Any element $br\in P^g(3,3), b\in B(2,3)$ can be represented by a simple curve in its
conjugation class unless $b$ is in the center of $B(2,3)$.

\end{lemma}
We have $ br = bx_1x_2$. The element $bx_1$ is represented by a simple curve $B$. Let $S$ be curve which contains $I_1$ and intersects $B$ transversally at exactly one point inside $I_1$.
The complimentary $ T - ( S- I_1)$ defines another group $B(3,2)$.
We assume that $p_0$ is not in $S$ and hence $bx_1$ is equivalent
 to a simple curve  with a desired orientation with respect to $x_2$.
That means $ (bx_1)x_2$ can be realized by a class of simple curve
 in $ P^g(3,3)$.
The classes $x_i$ are also realized by simple curves .   

$\Box$

\begin{corr} The elements $(br)^3 = 1$ for any $b\neq c = x_1x_2$.
\end{corr}

Indeed if $b$ is not in the center then $br$ is realized by a simple curve and we get the  result. The element $c^{-1}r = x_2^{-1}x^{-1}x_1 x_2^{-1} = x_2^{-2} = x_2$ in $P^g(3,3)$.

In dealing with the extension of $B(2,3)$ by an element $r$ we
will be using the following general  argument. Let $G$ be a finite group of exponent $3$ and $G'$ is obtained
from $G$ by adding $r$
and some relations of type $(br)^3 = 1$. Then the kernel of
 a natural projection $p: G'\to G$ is generated by the elements
$r^a = ara^{-1}, a \in G$. The group $G$ acts on this set
 of elements by left  translation $ g : r^a \to r^{ga}$. 
Any relation $(br)^3 = 1$ implies the relation :
$1 = brbr br = brb^{-1} b^{-1}r b r = r^b r^{b^{-1}}r =1$ and
similar relation for left translations of the orbits of
the cyclic group $B = (1,b,b^{-1})$. If in addition 
$(b^{-1}r)^3 = 1$ all the elements $r^b,r^{b^{-1}}, r$ commute
and any pair of them generate the same abelian group.

\begin{lemma} The group $P^g(3,3)$ is an abelian extension
of $B(2,3)$.
\end{lemma}

{\bf Proof} The kernel of the projection $p_r :P^g(3,3)\to B(2,3)$
which maps $r$ to $1$
is generated by the elements $r^b =  brb^{-1}$. All these  elements
commute with $r$ unless $b = c$. We also have
 $r^b r r^{b^{-1}} = 1$  and they all commute if $b$ is not in the center. Therefore $ r^a$ commutes with $r$ if it commutes with
 both $r^b,r^{b^{-1}}$. On the other hand $ r^{ab},r^{a}, r^{ab^{-1}}$
 commute if $r^b,r, r^{b^{-1}}$ commute.
Let $a$ be a generator of $B(2,3)$. Then $ r^{ca } r^{c} r^{ca^{-1}}= 1$  and all these elements
 commute, but $r^{ca}, r^{ca^{-1}}$ commute with $r$. This  implies
that all the elements $r^b, b\in B(2,3)$ commute with $r$.
After translation by  $B(2,3)$ we obtain that all
 the elements $r^g$ commute.

$\Box$

\begin{lemma}Let $T$ be the group generated by $r^b, b\in B(2,3)$.
Then $T$ is an  abelian group with $4$
generators.
\end{lemma}

{\bf Proof} Indeed the set of cyclic subgroups which don't lie in
the center generate a family of relations. Since we have established that
$T$ is an abelian group we shall write them in the additive form
$ r^{x}+r^{ax}+ r^{a^2 x} = 0$, $x \in B(2,3) a$ but not in the center $C$.
Denote by $T_S$ a subgroup of $T$ generated by a subset
 $S\subset B(2,3)$.
Let $A$ be  an abelian subgroup   generated by $a$ and $c$ which generates the center of $B(2,3)$. 
The summation over orbits of cyclic
 noncentral subgroups 
 gives zero.
Thus if we consider $T_A$ modulo a subgroup $T_C$ corresponding
 to the center we obtain $r^g + r^{gh} = 0, g\notin C$.
Hence the elements $r^g = -r^{g^{-1}}$ and $r^{gc}= r^{g}$ modulo
subgroup $r^c, r$. We obtain that $r^{a}$ generates $T_A$ modulo the  subgroup
$(r_c,r) $ and a sum over any orbit of $C$ is also zero.

Thus we have 
 $ r(x^{-1}) = - r - r^{x}$ 
and $ r^{x} + r^{y} +  r^{y^{-1}x^{-1}} = 0$ for the elements
 $x,y$ which lie in one abelian subgroup of $B(3,2)$.
The same is true modulo $r^c$ since we can apply the same
 argument to the quotient of $B(3,2)$ by the center.
Hence $r^{x} + r^{y} - r = r^{xy}$ (modulo($r^c$))
 for any $x,y\in B(2,3) $. 
In particular $ r^{a}, r^{b}, r^{c},r$ generate the group $T$.

$\Box$

\begin{lemma} $T$ is an elementary abelian group.
\end{lemma}

{\bf Proof} We have $r^{x^2} = 2 r^x - r$ and $r = 3 r^x- 2r$.
Hence $3r = 3r^x$ for any $x$.
Hence $ 3(r - r^x) = 0$ for any $x$.
Thus the elements of zero degree in $T$ constitute an
 elementary $3$-group $T_0$ which is a normal subgroup 
of $P^g(3,3)$. The quotient $T/T_0$ is a cyclic
group. The group 
  $P^g(3,3)/T_0$ is a central extension of
  $B(2,3)$. Since  $P^g(3,3)$ is generated by
 elements of order $3$ we obtain that $T/T_0 =\Bbb{Z}_3$.
  
$\Box$

\begin{corr} The number of elements in $T$ is $3^4$.
\end{corr}

 Hence the number of elements in $P^g(3,3)$ is equal to  $ 3^7$  and coincides with the number of elements in $B(3,3)$. Since there exists a
 surjective map $p: P^g(3,3)\to B(3,3)$ the groups coincide. Thus we have  proved  that group $P^g(n,3)$  coincide with $B(n,3)$ for all $n$.

\bigskip

\noindent
F. Bogomolov, Courant Institute, NYU and Steklov Institute, Moscow. \\
bogomolo@MATH8.CIMS.NYU.EDU

\bigskip

\noindent
L. Katzarkov, UC Irvine. \\ lkatzark@math.uci.edu


\begin{thebibliography}{40}

\bibitem{ART}{\bf M. Artin} {\em Algebraic approximation of structures
over complete local rings.\/} Publ. Math. I.H.E.S. {\bf 36} (1969),
pp. 23-58.
\bibitem{BKS} {\bf F. Bogomolov, L. Katzarkov} {\em Symplectic fourfolds and projective surfaces.\/} to appear in Proceedings of Georgia Topology
conference, 1996. 

\bibitem{BH} {\bf F. Bogomolov, D. Husemoller }{\em Geometric properties of curves defined over number fields.\/}  University of Maryland,
preprint 1993

\bibitem{CM}{\bf F. Campana} {\em Remarques sur le groupes de
K\"{a}hleriennes nilpotent.\/} Ann. Scien. L'Ecole Norm. Sup. 28,
fasc.3, 1995, pp. 307-316.
\bibitem{FL}{\bf E. Ballico, F. Catanese, C. Ciliberto (Eds.)} {\em
Trento examples}, Classification of higher dimensional varieties,
Proceedings, Trento 1990, Springer LNM, 1515, 1992, pp.  134-139.
\bibitem{DK}{\bf H. Clemens} {\em Degeneration of K\"{a}hler
varieties.\/} Duke Mathematical Journal, v. 44, (1977),no. 2, pp. 215 -290.

\bibitem{Co}{\bf R. Coleman }{\em Ramified torsion points on curves.\/} Duke Mathemtical Journal, v.54 (1987), pp. 615 - 640

\bibitem{DM}{\bf P. Deligne,  D. Mumford} {\em Irreducibility of the
moduli space of curves,\/} Publ. Math. I.H.E.S. {\bf 6} 1969, pp.  75-109.
\bibitem{GOMPF}{\bf R. Gompf } {\em A new construction of symplectic
manifolds.\/} Annals of Math., 142, 1995, pp. 527-595.
\bibitem{LR}{\bf B. Lasell, M. Ramachandran} {\em Observations on
harmonic maps and singular varieties.\/} Ann. Sci. \'{E}cole
Norm. Sup. {\bf 29}, 1996, pp. 135-148.
\bibitem{LM}{\bf L. Katzarkov } {\em On the Shafarevich maps.\/} to
appear in Proc. Symp. in Pure Math., Proceedings of AMS Alg. Geometry
conference Santa Cruz 1995.
\bibitem{KR}{\bf L. Katzarkov, M. Ramachandran} {\em On the universal coverings of algebraic surfaces.\/} Preprint, 1995.
\bibitem{FT}{\bf L. Katzarkov, M. Ramachandran, T. Pantev} {\em
Geometric factorization of linear fundamental groups \/.}, in
preparation.
\bibitem{LN}{\bf L. Katzarkov } {\em Nilpotent groups and universal
coverings of smooth projective varieties.\/} Journal of Diff. Geom., vol. 45, 1997, pp. 336- 349. 
\bibitem{K1}{\bf J. Koll\'ar} {\em Shafarevich maps and plurigenera of
algebraic varieties.\/} Inventiones Matematicae, 113, Fasc. 1, 1993, no. 1,
pp.  177-215.
\bibitem{K2}{\bf J. Koll\'ar} {\em Shafarevich maps and automorphic
forms.\/} Princeton University Press, Princeton, NJ, (1995).
\bibitem{MS}{\bf W. Magnus, A. Karrass, D. Solitar} {\em Combinatorial group
theory. Presentations of groups in terms of generators and relations. \/}, second revised edition. Dover
Publications, Inc., New York, 1976. 
\bibitem{MUM}{\bf D. Mumford} {\em The topology of normal
singularities of an algebraic surface and a criterion for
simplicity.\/} Publ. Math. I.H.E.S.  {\bf 9} 1961, pp. 5-22.
\bibitem{N1}{\bf T. Napier} {\em Convexity properties of coverings of
smooth projective varieties.\/} Mathematische Annalen, 286, 1990,
pp. 433-479.
\bibitem{SIM}{\bf C. T. Simpson} {\em Constructing variations of Hodge
structures using Yang-Mills theory and applications to
uniformization.\/}, J. Amer. Math. Soc. 1 (1988), no. 4, pp. 867--918.
\bibitem{TOL}{\bf D. Toledo } {\em Projective surfaces with
nonresidually finite fundamental groups. \/}  Publ.
Math I.H.E.S., No. 77 (1993), pp. 103--119. 
\bibitem{YAU}{\bf S.T. Yau} {\em Uniformizations of geometric
structures.\/} Proc. of Symp. in Pure Mathematics, v. 48, (1988),
pp. 265-274.

\bibitem{Za}{\bf Yu. Zarhin}  {\em Endomorphisms and torsion of abelian variety.\/} Duke Mathematical Journal, v. 54 (1987), pp. 131-145.

\end{thebibliography}
\end{document}